\documentclass[12pt,reqno]{amsart}
\usepackage{amssymb}
\usepackage{latexsym}
\usepackage[dvips]{graphics}
\usepackage{epsfig}
\usepackage{amsmath,amsfonts,amsthm,amssymb,amscd}
\input amssym.def
\input amssym.tex

\newcommand{\ncr}[2]{{#1 \choose #2}}

\newcommand{\twocase}[5]{#1 \begin{cases} #2 & \text{{\rm #3}}\\ #4
&\text{{\rm #5}} \end{cases}   }

\newcommand{\lp}{\left(}
\newcommand{\rp}{\right)}

\newcommand{\mattwo}[4]
{\left(\begin{array}{cc}
                        #1  & #2   \\
                        #3 &  #4
                          \end{array}\right) }

\newcommand\be{\begin{equation}}
\newcommand\ee{\end{equation}}
\newcommand\bp{\begin{proof}}
\newcommand\ep{\end{proof}}

\newcommand\bea{\begin{eqnarray}}
\newcommand\eea{\end{eqnarray}}
\newcommand\bml{\begin{multline}}
\newcommand\eml{\end{multline}}
\newcommand\bal{\begin{align}}
\newcommand\eal{\end{align}}
\newcommand\bi{\begin{itemize}}
\newcommand\ei{\end{itemize}}
\newcommand\ben{\begin{enumerate}}
\newcommand\een{\end{enumerate}}
\newcommand\bc{\begin{center}}
\newcommand\ec{\end{center}}
\newcommand\ba{\begin{array}}
\newcommand\ea{\end{array}}

\newtheorem{thm}{Theorem}[section]

\newtheorem{defi}[thm]{Definition}

\newcommand{\R}{\ensuremath{\mathbb{R}}}
\newcommand{\N}{\ensuremath{\mathbb{N}}}
\newcommand{\E}{\ensuremath{\mathbb{E}}}

\newcommand{\C}{\ensuremath{\mathbb{C}}}
\newcommand{\Z}{\ensuremath{\mathbb{Z}}}
\newcommand{\Q}{\mathbb{Q}}

\newcommand{\gep}{\epsilon}   
\newcommand{\gl}{\lambda}

\newcommand{\hphi}{\widehat{\phi}}  

\newcommand{\foh}{\frac{1}{2}}  

\newcommand{\logpm}{ \frac{\log p}{\log(m/\pi)} }      
\newcommand{\logptm}{ 2\frac{\log p}{\log(m/\pi)} }    
\newcommand{\fommt}{\frac{1}{m-2}}                             

\numberwithin{equation}{section}

\begin{document}

\title[Nuclei, Primes and the Random Matrix Connection]{Nuclei, Primes and the Random Matrix Connection}

\author{Frank W. K. Firk}
\email{fwkfirk@aol.com} \address{The Henry Koerner Center for Emeritus Faculty,
Yale University, New Haven, CT 06520}

\author{Steven J. Miller}
\email{Steven.J.Miller@williams.edu} \address{Department of Mathematics and Statistics, Williams College,\newline Williamstown, MA 01267}


\date{\today}

\begin{abstract}

In this article, we discuss the remarkable connection between two very different fields,
number theory and nuclear physics.  We describe the essential aspects of these fields, the quantities studied, and how insights in one have been fruitfully applied in the other. The exciting branch of modern mathematics -- random matrix theory -- provides the connection between the two fields. We assume no detailed knowledge of  number theory, nuclear physics, or random matrix theory; all that is required is some familiarity with linear algebra and probability theory, as well as some results from complex analysis.  Our goal is to provide the inquisitive reader with a sound overview of the subjects, placing them in their historical context in a way that is not traditionally given in the popular and technical surveys. 

\end{abstract}

\thanks{We would like to thank all of our colleagues in nuclear physics, number theory and random matrix theory for years of fruitful discussions. We also thank Mike Buchanan, Dan File, John Goes and the referees for very careful readings of an earlier draft. The second named author was partially supported by NSF grants  DMS0600848 and DMS0850577.}

\maketitle

\tableofcontents


\section{Summary}\label{sec:summary}

In the early 1970's a remarkable connection was unexpectedly discovered between two very different fields, nuclear physics and number theory, when it was noticed that random matrix theory accurately modeled many problems in each. Random matrix theory was first used in the early 1900's in the study of the statistics of population characteristics \cite{Wis}. The field developed rapidly in the 1950's when it was found to describe the spacing distributions of adjacent resonances (of the same spin and parity) observed in the interaction of low energy neutrons with nuclei \cite{Wig5}, and it flourished in the 1970's following a chance encounter between Hugh Montgomery and Freeman Dyson \cite{Mon} (when they saw it also predicted answers to many of the most difficult problems in number theory).

In this review article we describe the subjects and the quantities studied, and how insights in one field have been fruitfully applied in the other. We assume no familiarity with either subject; for the most part, basic linear algebra and probability theory suffice (though we need some results from complex analysis on analytic continuation and contour integration for some of the number theory calcuations). As there are many mathematical surveys of the subject, as well as some popular accounts \cite{Ha,Roc} of how the connection between the fields was noticed, our goal is to explain the broad brushstrokes of the theory without getting bogged down in the technical details. For those interested in a more mathematical survey, we recommend \cite{Con2,Con3,FSV,KaSa2,KeSn3} (see also Section 1.8 of \cite{Meh2}). Our point is to give the flavor of the subject, and bring these amazing connections to the attention of a wide audience. We concentrate on a representative sample of results and problems, and urge the interested reader to sample the bibliography. In particular, though we discuss many of the current statistics studied, we discuss the computations in detail only for the classical problem of the density of normalized eigenvalues of real symmetric matrices and the $1$-level density for the family of Dirichlet $L$-functions. We chose these examples as the key steps in the analysis of these problems are similar to many others, but the mathematical prerequisites to follow the calculations are significantly less.

In particular, our choice means that there are many important topics which will have only the briefest of mention (if any). To do the field justice would require a significantly longer article than this. Our hope is that by keeping the pre-requisites modest a large audience will be able to appreciate the striking similarities between two very different fields, and get a sense as to the nature of the computations. There are a few places where real and complex analysis is used (Fourier transforms and the residue theorem), as well as some abstract algebra or group theory (mostly group homomorphisms from $(\Z/m\Z)^\ast$ to complex numbers of absolute value 1). We state all needed results, and when possible provide brief explanations and proofs. While in the last part of the paper we concentrate on Dirichlet $L$-functions, we do mention additional families of $L$-functions (the background material is more substantial here, and we give only the briefest mention of the needed facts).

The paper is organized as follows. In \S\ref{sec:intro} we first give some number theory preliminaries to set the stage, describing some of the problems researchers are interested in, and how they are connected with the zeros of the Riemannn zeta function. We mention the famous Riemann Hypothesis, and how its veracity is related to understanding the prime numbers. This provides the motivation for studying the behavior of these zeros. The amazing observation, first noticed in the 1970's, is that many properties of these zeros can be modeled by random matrix theory, which had enjoyed a remarkable success in modeling nuclear physics. We briefly describe random matrix theory and discuss why it is applicable to so many problems. In \S\ref{sec:nuclear} we describe some of the history of nuclear physics, concentrating on the experimental results which laid the groundwork for the introduction of random matrix theory. We then sketch the proof of one of the most important results in the subject, Wigner's semi-circle law, in \S\ref{sec:wigner}. While other results are more closely related to the number theory quantities we wish to study, we give this proof as it highlights in a very accessible manner the techniques needed to attack a variety of problems. We then return to number theory in \S\ref{sec:rmttonumb} and discuss some (but by no means all!) of the earliest applications of random matrix theory. We concentrate on the $1$-level density, highlighting the similarities between this calculation and the proof of Wigner's semi-circle law. We give an interpretation of our number theory results in the language of nuclear physics, and then conclude with a very brief summary of some of the current avenues being explored. \\

\emph{For the reader:} The core of the paper is Sections \ref{sec:intro} and \ref{sec:nuclear}, where we describe the number theory problems and the nuclear physics history which led to the development of random matrix theory, as well as briefly summarizing random matrix theory. Sections \ref{sec:wigner} and \ref{sec:rmttonumb} are more advanced (especially the latter), where we give details of the calculations. Many of the more technical comments and some proofs of claims are relegated to the footnotes; these may safely be skipped by the reader interested in the broad brushstrokes of the theory and subject. For the benefit of the reader, we have also included in the footnotes definitions and explanations of much of the assumed background material to help keep the paper accessible.


\section{Introduction}\label{sec:intro}

\subsection{Number Theory Preliminaries}\label{sec:ntprelims}

The primes\footnote{An integer $n \ge 2$ is prime if the only positive integers that divide it are 1 and itself; if $n \ge 2$ is not prime we say it is composite. The number 1 is neither prime nor composite, but instead is called a unit.} are the building blocks of number theory: every integer can be written uniquely as a product of prime powers \cite{HW}.\footnote{This property is called unique factorization, and is one of the most important properties of prime numbers. If 1 were considered a prime, then unique factorization would fail; for example, if 1 were prime then $3^5\cdot 7$ and $1^{2009} \cdot 3^5 \cdot 7$ would be two different factorizations of $1701$.} One of the most important questions we can ask about primes is also one of the most basic: how many primes are there at most $x$? In other words, how many building blocks are there up to a given point?

Euclid proved over 2000 years ago that there are infinitely many primes; so, if we let $\pi(x)$ denote the number of primes at most $x$, we know $\lim_{x\to\infty} \pi(x) = \infty$. Euclid's proof is still used in courses around the world.\footnote{The proof is by contradiction. Assume there are only finitely many primes. We label them $p_1$, $p_2$, $\dots$, $p_n$. Consider the number $m=p_1 p_2 \cdots p_n + 1$. Either this is a new prime not on our list, or it is composite. If it is composite, it must be divisible by some prime; however, it cannot be divisible by any prime on our list, as each of these give remainder 1. Thus $m$ is either prime or divisible by a prime not on our list. In either case, our list was incomplete, proving that there are infinitely many primes. With a little bit of effort, one can show that Euclid's argument proves there is a $C > 0$ such that $\pi(x) \ge C \log \log x$. The first few primes generated in this manner are 2, 3, 7, 43, 13, 53, 5, 6221671, 38709183810571, 139, 2801. A fascinating question is whether or not every prime is eventually listed (see \cite{Sl}).} Can we do better? How rapidly does $\pi(x)$ go to infinity? In particular, what can we say about $\lim_{x\to\infty} \pi(x)/x$, which represents the probability that a number at most $x$ is prime?

The answer is given by the Prime Number Theorem, which states the number of primes at most $x$ is ${\rm Li}(x) + o({\rm Li}(x))$, where ${\rm Li}(x) = \int_2^x dt/\log t$ and for $x$ large, ${\rm Li}(x)$ is approximately $x/\log x$.\footnote{The notation $f(x) = o(g(x))$ means $\lim_{x\to\infty} f(x)/g(x) = 0$.} While it is possible to prove the prime number theorem elementarily \cite{Erd,Sel2}, the most informative proofs use complex numbers\footnote{A complex number $z \in \C$ is of the form $z = x+iy$. We call $x$ the real part and $y$ the imaginary part; we frequently denote these by $\Re(z)$ and $\Im(z)$, respectively.} and complex analysis, and lead to the fascinating connection between number theory and nuclear physics. One of the most fruitful approaches to understanding the primes is to understand properties of the Riemann zeta function, $\zeta(s)$, which is defined for $\Re(s) > 1$ by \be \zeta(s) \ = \ \sum_{n=1}^\infty \frac1{n^s}; \ee the series converges for $\Re(s) > 1$ by the integral test.\footnote{Let $s = \sigma+it$. Then $|\zeta(s)| \le \sum_{n=1}^\infty n^{-\sigma}$. The integral test from calculus states that this series converges if and only if $\int_1^\infty x^{-\sigma} dx$ converges, and this integral converges if $\sigma > 1$.} By unique factorization, it turns out that we may also write $\zeta(s)$ as a product over primes;\footnote{To see this, use the geometric series formula (see Footnote \ref{tetetet}) to expand $(1-p^{-s})^{-1}$ as $\sum_{k=0}^\infty p^{-ks}$ and note that $n^{-s}$ occurs exactly once on each side (and clearly every term from expanding the product is of the form $n^{-s}$ for some $n$.} this is called the Euler product of $\zeta(s)$, and is one of its most important properties:\footnote{\label{footnote:importanceeulerprod} We give two quick proofs of the importance of the Euler product by showing how it implies there are infinitely many primes. The first is $\sum 1/n^s \to \infty$ as $s\to 1+$, which means there must be infinitely many primes as otherwise the product is finite. The second proof is to note $\sum 1/n^2 = \pi^2/6$ is irrational; if there were only finitely many primes than the product would be rational. See for example \cite{MT-B} for details.} \be \zeta(s) \ = \ \sum_{n=1}^\infty \frac1{n^s} \ = \ \prod_{p\ {\rm prime}} \left(1 - \frac1{p^s}\right)^{-1}. \ee Initially defined only for $\Re(s) > 1$, using complex analysis the Riemann zeta function can be meromorphically continued\footnote{\label{tetetet}The subject of meromorphic continuation belongs to complex analysis. For the benefit of the reader who hasn't seen this, we give a brief example that will be of use throughout this paper, namely the geometric series formula $1+r+r^2 + r^3 + \cdots = 1/(1-r)$. Note that while the sum makes sense only when $|r| < 1$, $1/(1-r)$ is well-defined for all $r\neq 1$ and agrees with the sum whenever $|r| < 1$. We say $1/(1-r)$ is a meromorphic continuation of the sum.}  to all of $\C$, having only a simple pole with residue 1 at $s=1$. It satisfies the functional equation\footnote{One proof is to use the Gamma function, $\Gamma(s) = \int_0^\infty e^{-t} t^{s-1}dt$.  A simple change of variables gives $\int_0^\infty x^{\foh s - 1} e^{-n^2 \pi x}dx$ $=$
$\Gamma\left(\frac{s}{2}\right) / n^s\pi^{s/2}$. Summing over $n$ represents a multiple of $\zeta(s)$ as an integral. After some algebra we find $\Gamma\left(\frac{s}2\right) \zeta (s)$ $=$ $\int_1^\infty  x^{\foh s -1} \omega(x)dx$ $+$ $\int_1^\infty x^{-\foh s -1}\omega\left(\frac{1}{x}\right)dx$, with $\omega(x) = \sum_{n=1}^{+\infty} e^{-n^2 \pi
x}$. Using Poisson summation, we find $\omega\left(\frac{1}{x}\right) = -\foh + -\foh x^\foh + x^\foh \omega(x)$, which yields $\pi^{-\foh s}
\Gamma\left(\frac{s}2\right) \zeta (s)$ $=$ $\frac{1}{s(s-1)}$ $+$
$\int_1^\infty  (x^{\foh s -1}+ x^{- \foh s - \foh}) \omega(x)dx$, from which the claimed functional equation follows.}  \cite{Da,MT-B} \be\label{eq:defixi} \xi(s) \ =\ \foh s(s-1)
\Gamma\left(\frac{s}{2}\right)\pi^{-\frac{s}{2}} \zeta(s) \ = \ \xi(1-s).\ee The distribution of the primes is a difficult problem; however, the distribution of the positive integers is not, and has been completely known for quite some time! The hope is that we can understand $\sum_n 1/n^s$ as this involves sums over the integers, and somehow pass this knowledge on to the primes through the Euler product (see Footnote \ref{footnote:importanceeulerprod} for two examples).

Riemann \cite{Ri} (see \cite{Ed} for an English translation) observed a fascinating connection between the zeros of $\zeta(s)$ and the error term in the prime number theorem. As this relation is the starting point for our story, we describe the details in some length in the next paragraph. This part is a bit more technical and relies on complex analysis. \emph{The reader may safely skip most of the next paragraph; the key piece for the rest of the paper is \eqref{eq:keybeautyformriemannprimes}, where we show how the primes are connected to the zeros of $\zeta(s)$ (the function $\Lambda(n)$ which appears is defined in \eqref{eq:defnlambda}. } \\

One of the most natural things to do to a complex function is to take contour integrals of its logarithmic derivative \cite{La,SS}; this will yield information about zeros and poles (we'll see later that we can get even more information if we weight the integral with a test function). There are two expressions for $\zeta(s)$; however, for the logarithmic derivative it is clear that we should use the Euler product over the sum expansion, as the logarithm of a product is the sum of the logarithms. Let  \be\label{eq:defnlambda} \twocase{\Lambda(n)\  = \ }{\log p}{if $n=p^r$ for some integer $r$}{0}{otherwise.} \ee We find \be\label{eq:logderiv} \frac{\zeta'(s)}{\zeta(s)} \ = \ -\sum_p \frac{\log p \cdot p^{-s}}{1-p^{-s}} \ = \ -\sum_{n=1}^\infty \frac{\Lambda(n)}{n^s}  \ee (this is proved by using the geometric series formula to write $(1-p^{-s})^{-1}$ as $\sum_{k=0}^\infty 1/p^s$, collecting terms and then using the definition of $\Lambda(n)$). Moving the negative sign over and multiplying by $x^s/s$, we find \be \frac1{2\pi i} \int_{(c)} -\frac{\zeta'(s)}{\zeta(s)} \frac{x^s}{s}\ ds \ = \ \frac1{2\pi i} \int_{(c)} \sum_{n \le x} \Lambda(n)\left(\frac{x}{n}\right)^s \frac{ds}{s}, \ee where we are integrating over some line $\Re(s) = c > 1$. The integral on the right hand side is $1$ if $n < x$ and $0$ if $n > x$ (by choosing $x$ non-integral, we don't need to worry about $x=n$), and thus gives $\sum_{n \le x} \Lambda(n)$. By shifting contours and keeping track of the poles and zeros of $\zeta(s)$, the residue theorem\footnote{Let $f$ be a meromorphic function with only finitely many poles on an open set $U$ which is bounded by a `nice' curve $\gamma$. Thus at each point $z_0 \in U$ we have $f(x) = \sum_{n=N}^\infty a_n (z-z_0)^n$ with $N > -\infty$. If $N > 0$ we say $f$ has a zero of order $N$. If $N < 0$ we say $f$ has a pole of order $-N$, and in this case we call $a_{-1}$ the residue of $f$ at $z_0$ (for clarity, we often denote this by ${\rm Res}(f,z_0)$). If $f$ does not have a pole at $z_0$, then the residue is zero. Our assumption implies that there are only finitely many points where the residue is non-zero. The residue theorem states $\frac1{2\pi i} \oint_\gamma f(z)dz = \sum_{z \in U} {\rm Res}(f,z)$. One useful variant is to apply this to $f'(z)/f(z)$, which then counts the number of zeros of $f$ minus the number of poles;  another is to look at $f'(z)/f(z) \cdot g(z)$ where $g(z)$ is analytic, which we will do later when stating the explicit formula for the zeros of the Riemann zeta function.} \cite{La,SS} implies that the left hand side is \be x - \sum_{\rho: \zeta(\rho) = 0} \frac{x^\rho}{\rho}; \ee the $x$ term comes from the pole of $\zeta(s)$ at $s=1$ (remember we count poles with a minus sign), while the $x^\rho/\rho$ term arises from zeros; in both cases we must multiply by the residue, which is $x^\rho/\rho$ (it can be shown that $\zeta(s)$ has neither a zero nor a pole at $s=0$).\footnote{Some care is required with this sum, as $\sum 1/|\rho|$ diverges. The solution involves pairing the contribution from $\rho$ with $\overline{\rho}$; see for example \cite{Da}.} The Riemann zeta function vanishes whenever $\rho$ is a negative even integer; we call these the \textbf{trivial} zeros. These terms contribute $\sum_{k=-1}^\infty x^{-2k}/(2k) = -\frac12 \log(1-x^{-2})$. This leads to the following beautiful formula: \be\label{eq:keybeautyformriemannprimes}  x - \sum_{\rho: \Re(\rho) \in (0,1)\atop \zeta(\rho)=0} \frac{x^\rho}{\rho} - \frac12\log(1-x^{-2}) \ = \ \sum_{n \le x} \Lambda(n)\ee If we write $n$ as $p^r$, the contribution from all $p^r$ pieces with $r \ge 2$ is bounded by $2x^{1/2}\log x$ for $x$ large,\footnote{To see this, note $\sum_{p^2 \le x} \log p  \le x^{1/2} \log x$, while the contribution from $n = p^r$ with $r \ge 3$ is bounded by $\sum_{r \ge 3} \sum_{p^r \le x} \log p  \le x^{1/3} 2 \log^2 x$ (this is because $p^r \le x$ implies $r \le \log_2 x \le 2 \log x$).} thus we really have a formula for the sum of the primes at most $x$, with the prime $p$ weighted by $\log p$. Through partial summation, knowing the weighted sum is equivalent to knowing the unweighted sum.\footnote{Partial summation is the discrete analogue of integration by parts \cite{MT-B}. In our case, $\sum_{p \le x} \log p \sim x$ is equivalent to $\sum_{p \le x} 1 \sim x/\log x$.}\\

We can now see the connection between the zeros of the Riemann zeta function and counting primes at most $x$. The contribution from the trivial zeros is well-understood, and is just $-\frac12\log(1-x^{-2})$. The remaining zeros, whose real parts are in $[0,1]$, are called the \textbf{non-trivial} or \textbf{critical} zeros. They are far more important and more mysterious. The smaller the real part of these zeros of $\zeta(s)$, the smaller the error. Due to the functional equation, however, if $\zeta(\rho) = 0$ for a critical zero $\rho$ then $\zeta(1-\rho) = 0$ as well.\footnote{Note this is only true for zeros in the critical strip, namely $0 \le \Re(\rho) \le 1$; for zeros outside the critical strip we can and do have zeros of $\zeta(s)$ not corresponding to zeros of $\zeta(1-s)$ because of poles of the Gamma function.} Thus the `smallest' the real part can be is 1/2. This is the celebrated \textbf{Riemann Hypothesis}.\footnote{The Riemann Hypothesis is probably the most important mathematical aside ever in a paper. Riemann \cite{Ed,Ri} wrote (translated into English; note when he talks about the roots being real, he's writing the roots as $1/2 + i\gamma$, and thus $\gamma \in \R$ is the Riemann Hypothesis): \emph{...and it is very probable that all roots are real. Certainly one would wish for a stricter proof here; I have meanwhile temporarily put aside the search for this after some fleeting futile attempts,
as it appears unnecessary for the next objective of my investigation.} Though not mentioned in the paper, Riemann had developed a terrific formula for computing the zeros of $\zeta(s)$, and had checked (but never reported!) that the first few were on the critical line $\Re(s) = 1/2$. His numerical computations were only discovered decades later when Siegel was looking through Riemann's papers.} It has a plethora of applications throughout number theory and mathematics; counting primes is but one of many.\footnote{The prime number theorem is in fact equivalent to the statement that $\Re(\rho) < 1$ for any zero of $\zeta(s)$. The prime number theorem was first proved independently by Hadamard \cite{Had} and de la Vall\'{e}e Poussin \cite{dlVP} in 1896. Each proof crucially used results from complex analysis, which is hardly surprising given that Riemann had shown $\pi(x)$ is related to the zeros of the meromorphic function $\zeta(s)$. It was not until almost 50 years later that Erd\"{o}s \cite{Erd} and Selberg \cite{Sel2} obtained elementary proofs of the prime number theorem (in other words, proofs that did not use complex analysis, which was quite surprising as the prime number theorem was known to be equivalent to a statement about zeros of a meromorphic function). See \cite{Gol2} for some commentary on the history of elementary proofs.} It is clear, however, that the distribution of the zeros of the Riemann zeta function will be of primary (in both senses of the word!) importance.

If we assume the Riemann Hypothesis, all the zeros in the critical strip ($0 \le \Re(\rho) \le 1$) lie on the critical line $\Re(s) = 1/2$, and it makes sense to talk about the distribution between adjacent zeros. The purpose of this note is to discuss one of the most powerful models used to predict the behavior of these zeros, namely random matrix theory. While other methods have since been developed, random matrix theory (which we describe in the next subsection) was the first to make truly accurate, testable predictions. The general idea is that the behavior of zeros of the Riemann zeta function are well-modeled by the behavior of eigenvalues of certain matrices. This idea had previously been successfully used to model the distribution of energy levels of heavy nuclei (some of the fundamental papers and books on the subject, ranging from experiments to theory, include \cite{BFFMPW,DLL,Dy1,Dy2,FLM,FRG,FKPT,Gau,HH,HPB,Hu,Meh1,Meh2,MG,Po,Wig1,Wig2,Wig3,Wig4,Wig5,Wig6}). We describe the development of random matrix theory in nuclear physics in detail in the next section, and then delve into more of the details of the connection between the two subjects in \S\ref{sec:wigner} and \S\ref{sec:rmttonumb}.

\subsection{Random Matrix Theory Preliminaries}\label{sec:rmtprelims}

Before describing what we mean by random matrix theory and random matrix ensembles (i.e., sets of matrices), we quickly review the needed analysis and probability material, and then in the next subsection discuss why random matrix theory is so applicable at modeling a variety of problems.

Let $p(x)$ be a continuous or discrete probability distribution. For notational convenience we assume $p$ is continuous and use integral notation below, though similar statements hold in the discrete case. This means $p(x) \ge 0$, $\int_{-\infty}^\infty p(x)dx = 1$, and if $X$ is a random variable with density $p$, then the probability $X$ takes on a value in $[a,b]$ is just $\int_a^b p(x)dx$.

Let $X$ be a random variable with density $p$. We define the $k$\textsuperscript{th} moment of $p$, denoted $\mu_k$ or $\E[X^k]$, by \be \mu_k \ = \ \E[X^k] \ := \ \int_{-\infty}^\infty x^k p(x) dx. \ee The zeroth moment is always 1, and the first moment is called the \emph{mean}. The second moment is related to the variance. Recall the variance $\sigma^2$ is defined by \be \sigma^2 \ = \ \E[(X-\mu)^2] \ = \ \int_{-\infty}^\infty (x-\mu)^2 p(x) dx, \ee and equals the second moment if the mean is zero. For convergence issues, we typically are interested in random variables with zero mean, variance 1 and finite higher moments.\footnote{By this we mean $\int_{-\infty}^\infty |x|^k p(x)dx < \infty$. If the integrand is not absolutely convergent, then the value could depend on how we tend to infinity. The standard example is the Cauchy distribution, $p(x) = (\pi (1 + x^2))^{-1}$. Note $\int_{-\infty}^\infty |x| p(x)dx = \infty$, $\lim_{A\to\infty} \int_{-A}^A x p(x)dx = 0$ and $\lim_{A\to\infty} \int_{-A}^{2A} x p(x)dx = \pi^{-1} \log 2$.} While at first it might seem restrictive to assume we have mean 0 and variance 1, this is actually equivalent to the first and second moments are finite.\footnote{The reason is we can always adjust our probability distribution, in this case, to have mean 0 and variance 1 by simple translations and rescaling. For example, if the density $p$ has mean $\mu$ and variance $\sigma$, then $g(x) = \sigma^{-1} p(\sigma x + \mu)$ has mean 0 and variance 1. Thus the third moment (or the fourth if the third vanishes) are the first moments that truly show the `shape' of the density.} The moments are extremely important for understanding a density. While it is not the case that the moments uniquely determine a probability distribution\footnote{For $x \in [0,\infty)$, consider $f_1(x) =
(2\pi x)^{-1/2} \exp\left(-(\log x)^2/2\right)$ and $f_2(x) = f_1(x)\left[1 + \sin(2\pi \log x)\right]$. For $r \in \N$, the $r$\textsuperscript{th} moment of $f_1$ and $f_2$
is $\exp(r^2/2)$. The reason for the non-uniqueness of moments is that
the moment generating function $M_f(t) = \int_{-\infty}^\infty \exp(tx) f(x)dx$ does not converge in a neighborhood of the origin. See \cite{CaBe}, Chapter 2.}, they do for sufficiently nice distributions. The situation is similar to the theory of Taylor series. It is sadly not the case that every `nice' function agrees with its Taylor series in an arbitrarily small neighborhood about the point of expansion, even if by `nice' we mean infinitely differentiable!\footnote{The standard example is the function $f(x) = \exp(-1/x^2)$ if $|x| > 0$ and $0$ otherwise. By using the definition of the derivative and L'Hopital's rule, we see $f^{(n)}(0) = 0$ for all $n$, but clearly this function is non-zero if $|x| > 0$. Thus the radius of convergence is zero! This example illustrates how much harder real analysis can be than complex analysis. There if a function of a complex variable has even one derivative then it has infinitely many, and is given by its Taylor series in some neighborhood of the point.} See \cite{ShTa,Si} for more details.

We can now describe random matrix theory and the ensembles we'll study. Consider a real symmetric matrix $A$, so \be A \ = \ {\left(\begin{array}{ccccc}
                        a_{11}  & a_{12} & a_{13} & \cdots & a_{1N}  \\
                        a_{12} &  a_{22} & a_{23} & \cdots & a_{2N} \\
                        \vdots &  \vdots & \vdots & \ddots & \vdots \\
                        a_{1N} & a_{2N} & a_{3N} & \cdots & a_{NN}
                          \end{array}\right) } \ = \ A^T, \ \ \ a_{ij}
                          = a_{ji}. \nonumber\ \ee

We fix a density $p$, and define \be \mbox{Prob}(A) \ = \ \prod_{1 \le i \le j \le N}
p(a_{ij}). \nonumber\ \ee This means \bea \text{Prob}\left(A: a_{ij} \in [\alpha_{ij},
\beta_{ij}]\right) \ = \ \prod_{1 \le i \le j \le N}
\int_{x_{ij}=\alpha_{ij}}^{\beta_{ij}} p(x_{ij}) dx_{ij}.
\nonumber\ \eea The goal is to understand properties of the eigenvalues of $A$. We accomplish this by studying a related measure where we place point masses at the normalized eigenvalues. We use the Dirac delta funtionals $\delta(x - x_0)$, which is a unit point mass at $x_0$. This means $\int f(x) \delta(x-x_0)dx = f(x_0)$.\footnote{\label{footnote:diractdeltafunctional} We can consider the probability densities
$A_N(x)dx$, where $A$ is a probability density and $A_N(x) = N\cdot A(Nx)$. As $N \to \infty$
almost all the probability (mass) is concentrated in a narrower
and narrower band about the origin; we let $\delta(x)$ be the limit
with all the mass at one point. It is a discrete (as opposed to
continuous) probability measure, with infinite density but finite
mass. Note that $\delta(x-x_0)$ acts like a unit point mass;
however, instead of having its mass concentrated at the origin, it
is now concentrated at $x_0$.}

To each real symmetric matrix $A$, we attach a probability measure\footnote{As $A$ is real symmetric, the eigenvalues are real (see Footnote 19) and thus this measure is well defined.}
\bea \mu_{A,N}(x) & \ = \ & \frac{1}{N} \sum_{i=1}^N \delta\left( x
- \frac{\lambda_i(A)}{2\sqrt{N}} \right); \eea in \S\ref{sec:wignsemicircstatement} we'll see why we are normalizing the eigenvalues as we have done here. This measure counts the number of normalized eigenvalues in an interval: \bea  \int_a^b
\mu_{A,N}(x)dx & \ = \ & \frac{\#\left\{\gl_i:
\frac{\gl_i(A)}{2\sqrt{N}}\ \in \ [a,b]\right\}}{N}. \eea Using the definition of the Dirac delta functional, the $k$\textsuperscript{th} moment, which we denote $M_{A,N}(k)$, is readily computed: \bea\label{eq:kthmomentmank}
M_{A,N}(k) & \ = \ & \frac{\sum_{i=1}^N
\lambda_i(A)^k}{2^k N^{\frac{k}2+1}}. \eea While this is a nice, explicit formula for the $k$\textsuperscript{th} moment, it seems useless as we do not know the location of the eigenvalues of $A$; we will see in \S\ref{sec:wigner} that this is not the case at all.

There are many other ensembles of matrices worth studying. In addition to real symmetric matrices, complex Hermitian and symplectic are frequently studied.\footnote{These ensembles have behavior that is often described by a parameter $\beta$, which is 1 for real symmetric, 2 for complex Hermitian and 4 for symplectic matrices.} In this paper we concentrate on real symmetric matrices.

Random matrix theory models the behavior of a system by an appropriate set of matrices. Specifically, we calculate some quantity (say the probability two normalized eigenvalues are less than half the average spacing apart) for each matrix and then average over all matrices in our family. The hope, which is born out in many cases (ranging from number theory to nuclear physics to bus routes in Cuernevaca, Mexico\footnote{See \cite{BBDS, KrSe}. We quote from \cite{BBDS} who quote from \cite{KrSe}: There \emph{is no covering company responsible for organizing the city transport. Consequently, constraints such as a time table that represents external influence on the transport do not exist.
Moreover, each bus is the property of the driver. The drivers try to maximize their income
and hence the number of passengers they transport. This leads to competition among the
drivers and to their mutual interaction. It is known that without additive interaction the
probability distribution of the distances between subsequent buses is close to the Poisonian
distribution and can be described by the standard bus route model.... A Poisson-like distribution
implies, however, that the probability of close encounters of two buses is high (bus
clustering) which is in conflict with the effort of the drivers to maximize the number of
transported passengers and accordingly to maximize the distance to the preceding bus. In
order to avoid the unpleasant clustering effect the bus drivers in Cuernevaca engage people
who record the arrival times of buses at significant places. Arriving at a checkpoint, the
driver receives the information of when the previous bus passed that place. Knowing the
time interval to the preceding bus the driver tries to optimize the distance to it by either
slowing down or speeding up.} The papers go on to show the behavior is well-modeled by random matrix theory (specifically, ensembles of complex Hermitian matrices)!}), is that these system averages are close to the behavior of the system of interest. We describe this correspondence in greater detail below.

\subsection{Why Random Matrix Theory}

Why do random matrix models have a chance of giving useful answers to questions in nuclear physics and other subjects?  We consider one of the central problems of classical mechanics, namely the orbits in a solar system. It is possible to write down a closed form solution in the special case when there are just two point masses interacting through gravity.\footnote{Explicitly, given two points with masses $m_1$ and $m_2$ and initial velocities $\vec{v}_1$ and $\vec{v}_2$ and located at $\vec{r}_1$ and $\vec{r}_2$, we can describe how the system evolves in time given that gravity is the only force in play.} The three body problem, however, defies closed form solutions.\footnote{While there are known solutions for special arrangements of special masses, three bodies in general position is still open; see \cite{Wh} for more details.} From physical grounds we know of course that there is a solution; however, for our solar system we cannot analyze the solution well enough to determine whether or not billions of years from now Pluto will escape from the sun's influence!\footnote{Whether or not Pluto will regain planetary status is an entirely different question.}

As difficult as the above problem is, the situation is significantly worse when we try to understand the
behavior of heavy nuclei. Uranium, for instance, has over $200$ protons and neutrons in its nucleus, each subject to and contributing to complex forces. If we completely understood the theory of the nucleus, we could predict the energy levels; sadly, we are far from a complete understanding! As we'll see in the next section, physicists were able to gain some insights into the nuclear structure by shooting neutrons into the nucleus and analyzing the results; however, a complete understanding of the nucleus was, and still is, lacking.

\begin{figure}
\begin{center}
\scalebox{1}{\includegraphics{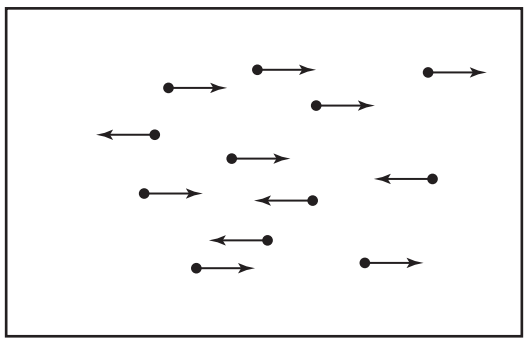}}
\caption{\label{fig:statmech}Molecules in a box}
\end{center}
\end{figure}

How should we attack such a problem? It's useful to recall other complex problems from physics and how they were successfully modeled. We consider a standard problem in statistical mechanics, namely calculating the pressure on a wall. Consider the box in Figure \ref{fig:statmech}.  For simplicity we assume that every molecule is moving either left or right, and all are traveling at the same speed. If we want to calculate the pressure on the left wall, we need to know how many particles strike the wall in an infinitesimal time. Thus  we need to know how many particles are close to the left wall and moving towards it. In a room there would be at least a mole (about $6.022 \cdot 10^{23}$) of air molecules, which means that this computation is well beyond our abilities.  Without going into all of the physics (see for example \cite{Re}), we can get a rough idea of what is happening. The complexity, the enormous number of configurations
of positions of the molecules, actually \emph{helps} us. For each configuration we can calculate the pressure due to that configuration. We then \emph{average} over all configurations, and it turns out that a generic configuration is close to the system average. This theory has enjoyed great success, and suggests a way to model nuclear physics.

Returning to our problem about heavy nuclei, from quantum mechanics we have the following equation governing our problem:
\bea H \Psi_n   \ = \  E_n \Psi_n, \eea where
$H$ is the Hamiltonian, whose entries depend on system,
$E_n$ are the energy levels and $\Psi_n$ are the energy eigenfunctions. Thus we have `reduced' nuclear physics to linear algebra. Unfortunately, there are two difficulties with this approach. The first is that $H$ is an infinite dimensional matrix, and the second is that we do not know any of the entries! This makes for quite a daunting task! Wigner's great insight was that this enormous complexity is similar to what we saw in Statistical Mechanics, and actually helps us. The interactions are so complex we might as well regard each entry as some randomly chosen number. Thus instead of considering the true $H$ for the system, we consider $N\times N$ real symmetric matrices with entries independently chosen from nice probability distributions. We compute whatever statistics we are interested in for these matrices, average over all matrices, and then take the $N\to\infty$ scaling limit. The main result is that
\emph{\textbf{the behavior of the eigenvalues of an arbitrary matrix
is often well approximated by the behavior obtained by averaging
over all matrices, and this is a good model for many systems, ranging from the energy levels of
heavy nuclei to the zeros of the Riemann zeta function.}}\footnote{This is reminiscent of the Central Limit Theorem. For example, if we average over all sequences of tossing a fair coin $2N$ times, we obtain $N$ heads, and \emph{most} sequences of $2N$ tosses will have approximately $N$ heads, where approximately means deviations on the order of $\sqrt{N}$.}


\section{Nuclear Physics History}\label{sec:nuclear}

Below we discuss some of the history of investigations of the nucleus, concentrating on the parts that led to the introduction of random matrix theory to the subject. We mention some of the connections with number theory, which will be explored in much greater detail later.

\subsection{Introduction}
	The Riemann Hypothesis asserts that the non-trivial zeros of the Riemann zeta function are of the form $\rho = 1/2+i\gamma_\rho$ with $\gamma_\rho$ real. About the year 1913, P\'{o}lya conjectured\footnote{The first reference to this conjecture in the literature might not have been until 1973 by Montgomery \cite{Mon}.}   that the $\gamma_\rho$ are the eigenvalues of a naturally occurring, unbounded, self-adjoint operator, and are therefore real.\footnote{If $v$ is an eigenvector with eigenvalue $\lambda$ of a Hermitian matrix $A$ (so $A = A^\ast$ with $A^\ast$ the complex conjugate transpose of $A$, then $v^\ast (Av) = v^\ast (A^\ast v) = (Av)^\ast v$; the first expression is $\lambda ||v||^2$ while the last is $\overline{\lambda} ||v||^2$, with $||v||^2 = v^\ast v = \sum |v_i|^2$ non-zero. Thus $\lambda = \overline{\lambda}$, and the eigenvalues are real. This is one of the most important properties of Hermitian matrices, as it allows us to order the eigenvalues.}\label{footnote:hermevaluesreal} Later, Hilbert contributed to the conjecture, and reportedly introduced the phrase `spectrum' to describe the eigenvalues of an equivalent Hermitian operator, apparently by analogy with the optical spectra observed in atoms.  This remarkable analogy pre-dated Heisenberg's Matrix Mechanics and the Hamiltonian formulation of Quantum Mechanics by more than a decade.  Not surprisingly, the P\'{o}lya-Hilbert conjecture was considered so intractable that it was not pursued for decades, and Random Matrix Theory remained in a dormant state.  To quote Diaconis \cite{Di1}: ``\emph{Historically, Random Matrix Theory was started by Statisticians \cite{Wis} studying the correlations between different features of population (height, weight, income...).  This led to correlation matrices with $(i, j)$ entry the correlation between the $i$th  and $j$th features.  If the data were based on a random sample from a larger population, these correlation matrices are random; the study of how the eigenvalues of such samples fluctuate was one of the first great accomplishments of Random Matrix Theory.}''   Diaconis \cite{Di2} has given an extensive review of Random Matrix Theory from the perspective of a statistician.  A strong argument can be made, however, that Random Matrix Theory, as we know it today in the Physical Sciences, began in a formal mathematical sense with the Wigner surmise \cite{Wig5} concerning the spacing distribution of adjacent resonances (of the same spin and parity) in the interactions between low-energy neutrons and nuclei, discussed below.

\subsection{Nuclear Physics and Random Matrix Theory}
   	
	The period from the mid-1930's to the late 1970's was the Golden Age of Neutron Physics; widespread interest in understanding the physics of the nucleus, coupled with the need for accurate data in the design of nuclear reactors, made the field of Neutron Physics of global importance in fundamental Physics, Technology, Economics, and Politics.
	In the mid-1950's, a discovery was made that turned out to have far-reaching consequences beyond anything that those working in the field could have imagined.  For the first time, it was possible to study the microstructure of the continuum in a strongly-coupled, many-body system, at very high excitation energies. This unique situation came about as the result of the following facts:\\

\bi

\begin{figure}
\begin{center}
\scalebox{.9}{\includegraphics{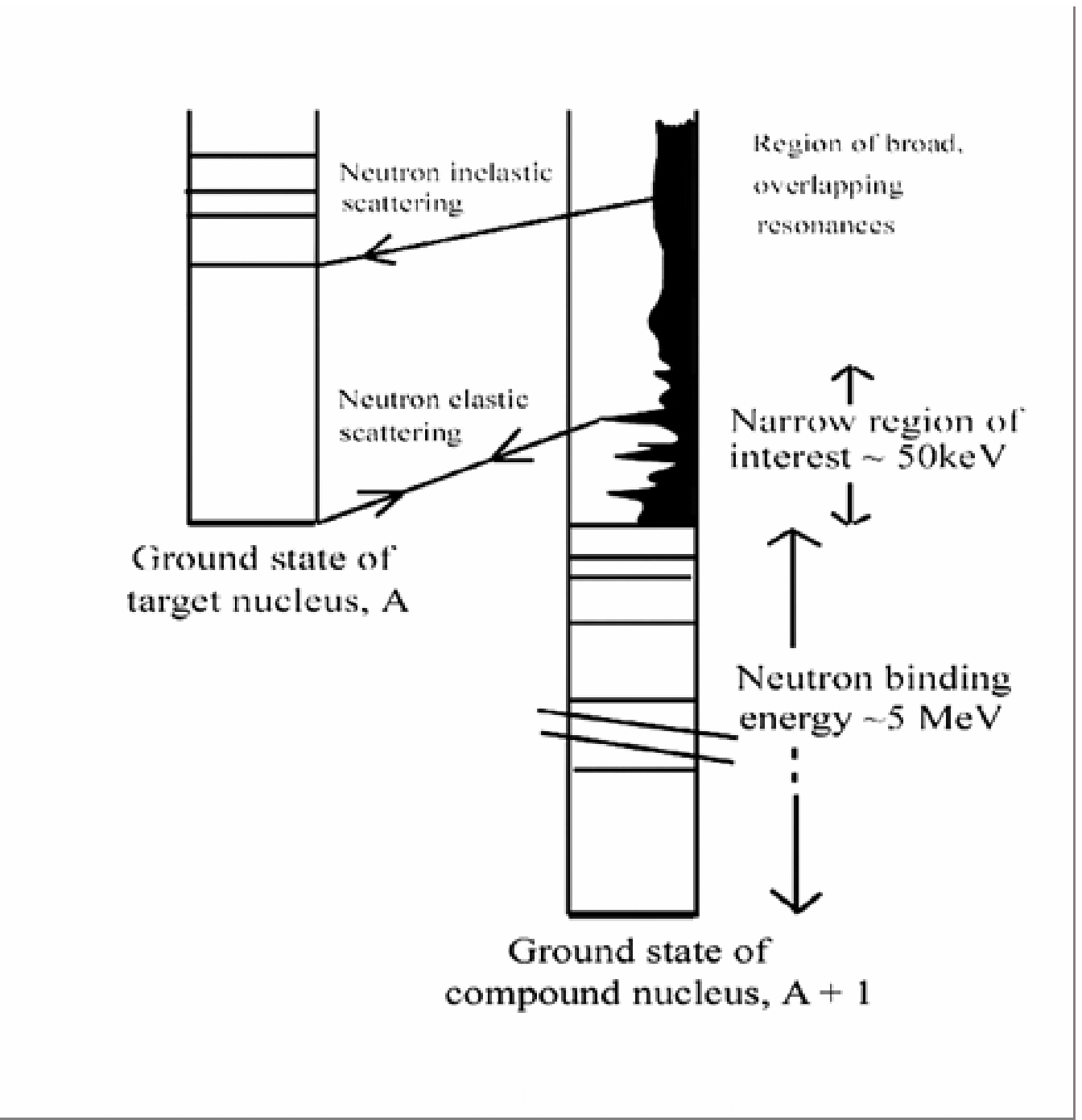}}
\caption{\label{fig:groundstatesneutrons} An energy-level diagram showing the location of  highly-excited resonances in the compound nucleus formed by the interaction of a neutron, n, with a nucleus of mass number A.  Nature provides us with a narrow energy region in which the resonances are clearly separated, and are observable.}
\end{center}
\end{figure}

\item	Neutrons, with kinetic energies of a few electron-volts, excite states in
compound nuclei at energies ranging from about 5 million electron-volts to almost 10 million electron-volts -- typical neutron binding energies. Schematically, see Figure \ref{fig:groundstatesneutrons}.

\ \\

\item	Low-energy resonant states in heavy nuclei (mass numbers greater than about 100) have lifetimes in the range $10^{-14}$ to $10^{-15}$ seconds, and therefore they have widths of about 1 eV.  The compound nucleus loses all memory of the way in which it is formed.  It takes a relatively long time for sufficient energy to reside in a neutron before being emitted.  This is a highly complex, statistical process.  In heavy nuclei, the average spacing of adjacent resonances is typically in the range from a few eV to several hundred eV. \\

\item	Just above the neutron binding energy, the angular momentum barrier restricts the possible range of values of total spin of a resonance, \textbf{J} (\textbf{J} = \textbf{I} + \textbf{i} + \textbf{l}, where \textbf{I} is the spin of the target nucleus, \textbf{i} is the neutron spin, and \textbf{l} is the relative orbital angular momentum).  This is an important technical point.	 \\

\item	The neutron time-of-flight method provides excellent energy resolution at energies up to several keV.  (See Firk \cite{Fi} for a review of time-of-flight spectrometers.) \\

\ei

 	A 1-eV neutron travels 1 meter in 72.3 microseconds.  At non-relativistic energies, the energy resolution $\Delta E$ at an energy $E$ is simply: \be \Delta E \approx 2E \Delta t/t_E, \ee
where $\Delta t$ is the \emph{total} timing uncertainty, and $t_E$ is the flight time for a neutron of energy $E$.

	In 1958, the two highest-resolution neutron spectrometers in the world had total timing uncertainties  $\Delta t \approx 200$ nanoseconds.  For a flight-path length of 50 meters the resolution was $\Delta E \approx 3$ eV at 1 keV.

	In ${\ }^{238}{\rm U} + {\rm n}$, the excitation energy is about 5 MeV; the effective resolution for a 1 keV-neutron was therefore \be
   \Delta E/E_{{\rm effective}} \ \approx \  6 \cdot 10^{-7}. \ee
(at 1 eV, the effective resolution was about $10^{-11}$).

	Two basic broadening effects limit the sensitivity of the method, they are:
\ben

\item Doppler broadening of the resonance profile due to the thermal motion of the target nuclei; it is characterized by the quantity $\delta \approx 0.3\sqrt{E/A}$  (eV), where $A$ is the mass number of the target.  If $E = 1$ keV and $A = 200$, $\delta \approx 0.7$ eV, a value that may be ten times greater than the natural width of the resonance.  \\

\item  Resolution broadening of the observed profile due to the finite resolving power of the spectrometer.  For a review of the experimental methods used to measure neutron total cross sections see Firk and Melkonian \cite{FM}.  Lynn \cite{Ly} has given a detailed account of the theory of neutron resonance reactions. \\

\een

	In the early 1950's, the field of low-energy neutron resonance spectroscopy was dominated by research groups working at nuclear reactors.  They were located at National Laboratories in the United States, the United Kingdom, Canada, and the former USSR.  The energy spectrum of fission neutrons produced in a reactor is moderated in a hydrogenous material to generate an enhanced flux of low-energy neutrons.  To carry out neutron time-of-flight spectroscopy, the continuous flux from the reactor is ``chopped'' using a massive steel rotor with fine slits through it.  At the maximum attainable speed of rotation (about 20,000 rpm), and with slits a few thousandths-of-an-inch in width, it is possible to produce pulses each with a duration approximately 1 $\mu$sec.  The chopped beams have rather low fluxes, and therefore the flight paths are limited in length to less than 50 meters.  The resolution at 1keV is then $\Delta E \approx 20$ eV, clearly not adequate for the study of resonance spacings about 10 eV.

	In 1952, there were only four accelerator-based, low-energy neutron spectrometers operating in the world.  They were at Columbia University in New York City, Brookhaven National Laboratory, the Atomic Energy Research Establishment, Harwell, England, and at Yale University. The performances of these early accelerator-based spectrometers were comparable with those achieved at the reactor-based facilities.  It was clear that the basic limitations of the neutron-chopper spectrometers had been reached, and therefore future developments in the field would require improvements in accelerator-based systems.

	In 1956, a new high-powered injector for the electron gun of the Harwell electron linear accelerator was installed to provide electron pulses with very short durations (typically less than 200 nanoseconds) \cite{FRG}.  The pulsed neutron flux (generated by the ($\gamma$, n) reaction) was sufficient to permit the use of a 56-meter flight path; an energy resolution of 3 eV at 1 keV was achieved.

	At the same time, Professors Havens and Rainwater (pioneers in the field of neutron time-of-flight spectroscopy) and their colleagues at Columbia University were building a new 385-MeV proton synchrocyclotron a few miles north of the campus (at the Nevis Laboratory). The accelerator was designed to carry out experiments in meson physics and low-energy neutron physics (neutrons generated by the (p, n) reaction).  By 1958, they had produced a pulsed proton beam with duration of 25 nanoseconds, and had built a 37-meter flight path \cite{RDRH, DRRH}.  The hydrogenous neutron moderator generated an effective pulse width of about 200 nanoseconds for 1 keV-neutrons.  In 1960, the length of the flight path was increased to 200 meters, thereby setting a new standard in neutron time-of-flight spectroscopy \cite{GRPH}.

\subsection{The Wigner Surmise}
	At a conference on Neutron Physics by Time-of-Flight, held in Gatlinburg, Tennessee on November 1st and 2nd, 1956, Professor Eugene Wigner (Nobel Laureate in Physics, 1963) presented his surmise regarding the theoretical form of the spacing distribution of adjacent neutron resonances (of the same spin and parity) in heavy nuclei.  At the time, the prevailing wisdom was that the spacing distribution had a Poisson form (see, however, \cite{GP}).  The limited experimental data then available was not sufficiently precise to fix the form of the distribution (see \cite{Hu}).  The following quotation, taken from Wigner's presentation at the conference, introduces the concept of random matrices in Physics, for the first time:

``\emph{Perhaps I am now too courageous when I try to guess the distribution of the distances between successive levels.  I should re-emphasize that levels that have different $J$-values (total spin) are not connected with each other.  They are entirely independent.  So far, experimental data are available only on even-even elements.  Theoretically, the situation is quite simple if one attacks the problem in a simple-minded fashion.  The question is simply `what are the distances of the characteristic values of a symmetric matrix with random coefficients?'}

\emph{We know that the chance that two such energy levels coincide is infinitely unlikely. We consider a two-dimensional matrix, $\mattwo{a_{11}}{a_{12}}{a_{21}}{a_{22}}$,
in which case the distance between two levels is $\sqrt{(a_{11} - a_{22})^2 + 4a_{12}^2}$.  This distance can be zero only if $a_{11} = a_{22}$ and $a_{12} = 0$. The difference between the two energy levels is the distance of a point from the origin, the two coordinates of which are $(a_{11} - a_{22})$ and $a_{12}$.  The probability that this distance is $S$ is, for small values of $S$, always proportional to $S$ itself because the volume element of the plane in polar coordinates contains the radius as a factor….}

\emph{The probability of finding the next level at a distance $S$ now becomes proportional to $SdS$. Hence the simplest assumption will give the probability}
\be                            \frac{\pi}{2} \rho^2 \exp\left( -\frac{\pi}{4}\rho^2S^2\right) SdS \ee
\emph{for a spacing between $S$ and $S + dS$.}

\emph{If we put $x = \rho S = S/\langle S \rangle$, where $\langle S \rangle$ is the mean spacing, then the probability distribution takes the standard form}
\be p(x)dx\  = \  \frac{\pi}{2}\ x\ \exp\left(-\pi x^2/4\right)dx, \ee
\emph{where the coefficients are obtained by normalizing both the area and the mean to unity.''  }

	This form, in which the probability of zero spacing is zero, is strikingly different from the Poisson form				 \be p(x)dx\ =\ \exp(-x)dx \ee
in which the probability is a maximum for zero spacing. The form of the Wigner surmise had been previously discussed by Wigner himself  \cite{Wig1}, and by Landau and Smorodinsky \cite{LS}, but not in the spirit of Random Matrix Theory.

	It is interesting to note that the Wigner distribution is a special case of a general statistical distribution, named after Professor E. H. Waloddi Weibull (1887-1979), a Swedish engineer and statistician \cite{Wei}.  For many years, the distribution has been in widespread use in statistical analyses in industries such as aerospace, automotive, electric power, nuclear power, communications, and life insurance.\footnote{In fact, one of the authors has used Weibull distributions to model run production in major league baseball, giving a theoretical justification for Bill James' Pythagorean Won-Loss formula \cite{Mil3}.}  The distribution gives the lifetimes of objects and is therefore invaluable in studies of the failure rates of objects under stress (including people!).  The Weibull probability density function is
\be {\rm Wei}(x; k, \lambda) = \frac{k}{\lambda} \left(\frac{x}{\lambda}\right)^{k-1} \exp\left(-(x/\lambda)^k\right) \ee
where $x \ge 0$, $k > 0$ is the \emph{shape} parameter, and $\lambda > 0$ is the \emph{scale} parameter.  We see that ${\rm Wei}(x; 2, 2/\sqrt{\pi}) = p(x)$, the Wigner distribution. Other important Weibull distributions are given in the following list
\bi
\item 	${\rm Wei} (x;1, 1) = \exp(-x)$ the Poisson distribution;

\item  ${\rm Wei} (x; 2, \lambda) = {\rm Ray}(\lambda)$, the Rayleigh distribution;

 \item ${\rm Wei} (x; 3, \lambda)$ is approximately a normal distribution.\footnote{Obviously this Weibull cannot be a normal distribution, as they have very different decay rates for large $x$, and this Weibull is a one-sided distribution! What we mean is that for $0 \le x \le 2$ this Weibull is well approximated by a normal distribution which shares its mean and variance, which are (respectively) $\Gamma(4/3) \approx .893$ and $\Gamma(5/3)-\Gamma(4/3)^2 \approx .105$.}

\ei

For ${\rm Wei}(x; k, \lambda)$, the mean is $\lambda \Gamma\left(1 + (1/k)\right)$,  the median is $\lambda \log(2)^{1/k}$,  and    the mode is $\lambda(k - 1)^{1/k}/k^{1/k}$, if $k>1$. As $k \to \infty$, the Weibull distribution has a sharp peak at $\lambda$.\footnote{Historically, Frechet introduced this distribution in 1927, and Nuclear Physicists often refer to the Weibull distribution as the Brody distribution \cite{BFFMPW}. }

	At the time of the Gatlinburg conference, no more than 20 s-wave neutron resonances had been clearly resolved in a single compound nucleus and therefore it was not possible to make a definitive test of the Wigner surmise.    Immediately following the conference, J. A. Harvey and D. J. Hughes \cite{HH}, and their collaborators, working at the fast-neutron-chopper-groups at the high flux reactor at the Brookhaven National Laboratory, and at the Oak Ridge National laboratory, gathered their own limited data, and all the data from neutron spectroscopy groups around the world, to obtain the first \emph{global spacing distribution} of s-wave neutron resonances.  Their combined results, published in 1958, showed a distinct lack of very closely spaced resonances, in agreement with the Wigner surmise.

	By late 1959, the experimental situation had improved, greatly.  At Columbia University, two students of Professors Havens and Rainwater completed their PhD theses; one, Joel Rosen \cite{RDRH}, studied the first 55 resonances in ${\ }^{238}{\rm U} + {\rm n}$ up to 1 keV, and the other, J Scott Desjardins \cite{DRRH}, studied resonances in two silver isotopes (of different spin) in the same energy region.  These were the first results from the new high-resolution neutron facility at the Nevis cyclotron.

\begin{figure}
\begin{center}
\scalebox{.9}{\includegraphics{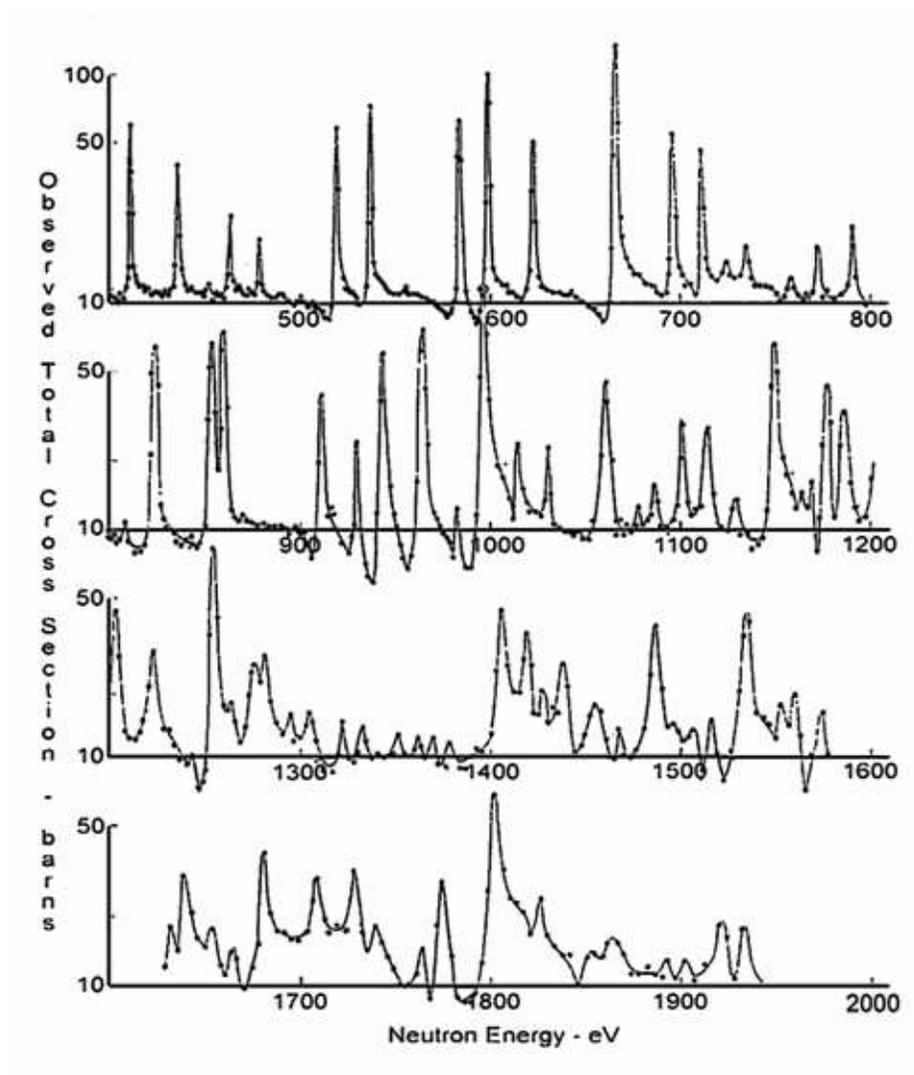}}
\caption{\label{fig:highresonancestudies} High resolution studies of the total neutron cross section of ${\ }^{238}{\rm U}$, in the energy range 400eV - 1800eV (12).  The vertical scale (in units of "barns") is a measure of the effective area of the target nucleus.}
\end{center}
\end{figure}

	At Harwell, Firk, Lynn, and Moxon \cite{FLM} completed their study of the first 100 resonances in ${\ }^{238}{\rm U} + {\rm n}$ at energies up to 1.8 keV; their measurement of the total neutron cross section for the interaction ${\ }^{238}{\rm U} + {\rm n}$  in the energy range 400--1800 eV is shown in Figure \ref{fig:highresonancestudies}.

When this experiment began in 1956, no resonances had been resolved at energies above 500 eV.  The distribution of adjacent spacings of the first 100 resonances in the single compound nucleus, ${\ }^{238}{\rm U} + {\rm n}$, ruled out an exponential distribution and provided the best evidence (then available) in support of Wigner's proposed distribution.	

\begin{figure}
\begin{center}
\scalebox{.9}{\includegraphics{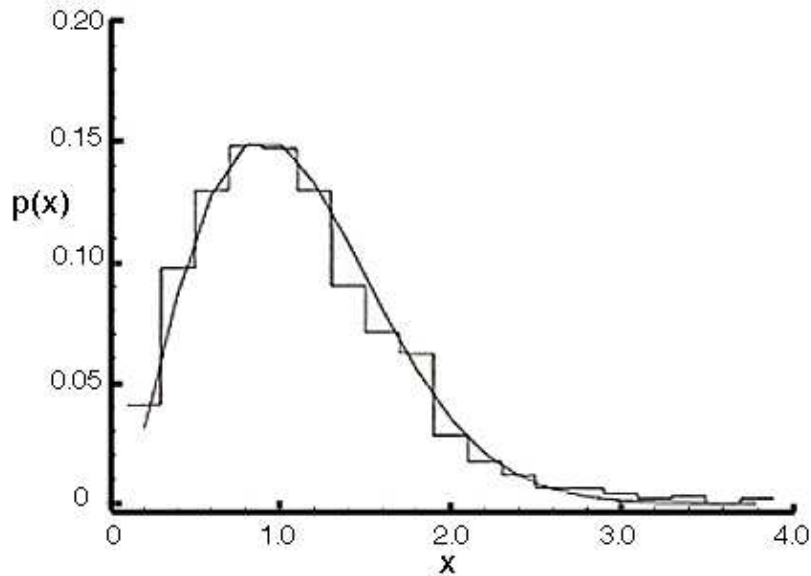}}
\caption{\label{fig:firkfig3} A Wigner distribution fitted to the spacing distribution of 932 s-wave resonances in the interaction ${\ }^{238}{\rm U} + {\rm n}$ at energies up to 20 keV.}
\end{center}
\end{figure}

	Over the last half-century, numerous studies have not changed the basic findings.  At the present time, almost 1000 s-wave neutron resonances in the compound nucleus ${\ }^{239}{\rm U}$ have been observed in the energy range up to 20 keV.  The latest results, with their greatly improved statistics, are shown in Figure \ref{fig:firkfig3} \cite{DLL}.	 

\subsection{Further Developments}

	The first numerical investigation of the distribution of successive eigenvalues associated with random matrices was carried out by Porter and Rozenzweig in the late 1950's \cite{PR}. They diagonalized a large number of matrices where the elements are generated randomly but constrained by a probability distribution.  The analytical theory developed in parallel with their work: Mehta \cite{Meh1}, Mehta and Gaudin \cite{MG}, and Gaudin \cite{Gau}. At the time it was clear that the spacing distribution was not influenced significantly by the chosen form of the probability distribution. Remarkably, the $n \times n$ distributions had forms given \emph{almost exactly} by the original Wigner $2\times 2$ distribution.

	The linear dependence of $p(x)$ on the normalized spacing $x$ (for small $x$) is a direct consequence of the \emph{symmetries} imposed on the (Hamiltonian) matrix, $H(h_{ij})$.  Dyson [Dy1] discussed the general mathematical properties associated with random matrices and made fundamental contributions to the theory by showing that different results are obtained when different \emph{symmetries} are assumed for $H$.  He introduced three basic distributions; in Physics, only two are important, they are:\\

\bi
\item	the Gaussian Othogonal Ensemble (GOE) for systems in which rotational symmetry and time-reversal invariance holds (the Wigner distribution): $p(x) = (\pi/2) x \exp\left(-(\pi/4)x^2\right)$;\\

\item the Gaussian Unitary Ensemble (GUE) for systems in which time-reversal invariance does not hold (French \emph{et al.} \cite{FKPT}):	$p(x) = (32/\pi^2) x^2\exp(-(\pi/4)x^2)$.\\

\ei

The mathematical details associated with these distributions are given in \cite{Meh1}.

 	The impact of these developments was not immediate in Nuclear Physics.  At the time, the main research endeavors were concerned with the structure of nuclei--experiments and theories connected with Shell-, Collective-, and Unified models, and with the nucleon-nucleon interaction.  The study of Quantum Statistical Mechanics was far removed from the main stream.   Almost two decades went by before Random Matrix Theory was introduced in other fields of Physics (see, for example, Bohigas, Giannoni and Schmit \cite{BGS} and Alhassid \cite{Al}).

\subsection{From Physics to Number Theory}

Interestingly, the next development occurred in an area having nothing to do with Physics.  In the field of Number Theory, perhaps the greatest unsolved problem has to do with the Riemann conjecture (that dates from the mid-19\textsuperscript{th} century): if $\zeta(s) = \sum 1/n^s$, then every complex number $\rho$ in the critical strip ($0 \le \Re(\rho) \le 1$) at which the analytic continuation of $\zeta(s)$ has a non-trivial zero has real part equal to $1/2$.  In 1914, Hardy \cite{Har} proved that there are infinitely many zeros of the zeta function on the line critical line $\Re(s) = 1/2$. Later Selberg \cite{Sel1} proved a small positive percentage are on this line; this was improved by Levinson \cite{Lev} to a third, and now thanks to Conrey \cite{Con1} we know at least two-fifths lie on the line.\footnote{There is an interesting perspective to proving more than a third of the zeros lie on the critical line. As zeros off the line occur in complex conjugate pairs, proving more than a third of all non-trivial zeros lie on the line is equivalent to more than a half of all zeros with real part at least 1/2 are on the line! Thus, in this sense, a `majority' of all zeros are on the critical line.}

In the early 1970's, Hugh Montgomery, a mathematician at the University of Michigan, was investigating the relative spacing of the zeros of the zeta function \cite{Mon} (because of applications to the class number problem\footnote{The class number measures how much unique factorization fails in the ring of integers of a finite extension of $\Q$, and thus is an extremely important property of these fields. For example, $\Z[i] = \{a+ib: a, b \in \Z\}$ has unique factorization, while $\Z[i\sqrt{5}] = \{a+ib\sqrt{5}: a, b \in \Z\}$ does not (in the latter, note we can write 6 as either $2 \cdot 3$ or $(1+i\sqrt{5})(1-i\sqrt{5})$, and none of these four numbers can be factored as $(a+ib\sqrt{5})(c+id\sqrt{5})$ without one of the two factors being a unit (the units are numbers in the ring whose norm is 1; in $\Z[i\sqrt{5}]$, these numbers are $\pm 1$ (in $\Z[\sqrt{5}]$ we would also have numbers such as $2+\sqrt{5}$, as $(2+\sqrt{5})(-2+\sqrt{5}) = 1$). The class number problem is to find all imaginary quadratic fields with a given class number; see \cite{St,Wa} for more details and results. It turns out (see \cite{CI}) that if there are many small (relative to the average spacing) gaps between zeros of $\zeta(s)$ on the critical line, then there are terrific lower bounds for the class number; another connection between the class number and zeros of $L$-functions is through the work of Goldfeld \cite{Gol1} and Gross-Zagier \cite{GZ}.}).  Let us recall that, if we have a series of points distributed randomly along a line, with average density normalized to 1, and we treat the coordinates of the points as independent random variables, then the probability of finding $j$ points in a given interval of length $x$ is the Poisson distribution \be	 \frac{x^j}{j!} \exp(-x). \ee

For our real symmetric and complex Hermitian random matrix ensembles, the probability of finding more than one eigenvalue in a short interval is less than that given by the Poisson distribution -- the eigenvalues of the random matrix are said to `repel' each other.  The pair and higher level correlation function describe this effect (we discuss these functions in greater detail later in the paper; knowing all the correlation functions is equivalent to knowing the neighbor spacings).  Montgomery studied the pair correlation function for the zeros of the zeta function and he gave evidence that it has the asymptotic form be    \be   1 - \left(\frac{\sin \pi x}{\pi x}\right)^2. \ee

\begin{figure}
\begin{center}
\scalebox{.9}{\includegraphics{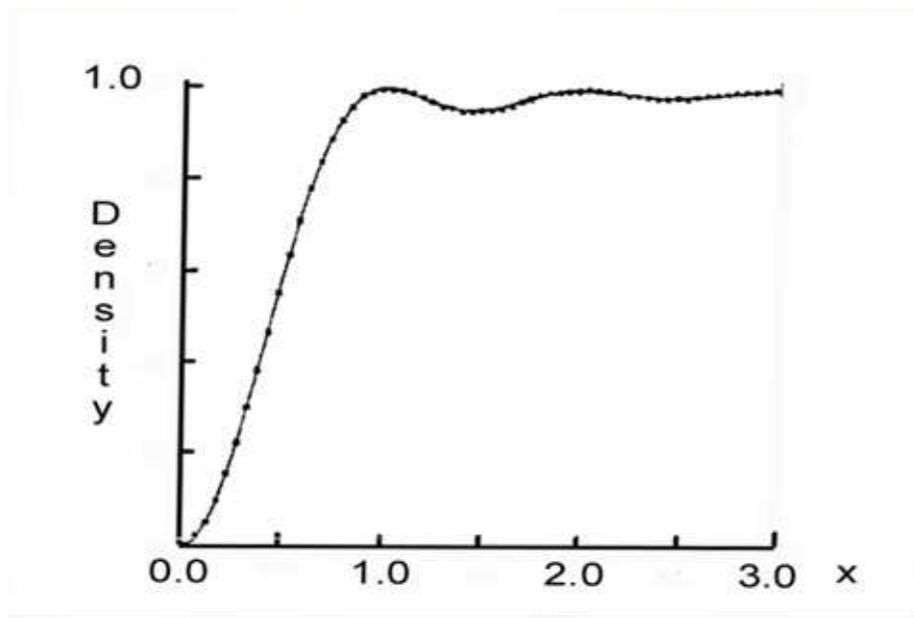}}
\caption{\label{fig:odlyzko} Odlyzko's test of the Montgomery conjecture,
involving 70 million Riemann zeros near $10^{20}$.}
\end{center}
\end{figure}

	At a chance meeting between Montgomery and Dyson at Princeton in the early 1970's, Montgomery showed his pair correlation function to Dyson, who recognized it as the pair correlation function of eigenvalues of random Hermitian matrices in a Gaussian Unitary Ensemble (an ensemble without time-reversal invariance).  In a masterful numerical calculation of the distribution of spacings between zeros of the zeta function, Andrew Odlyzko \cite{Od1, Od2} tested the Montgomery conjecture by studying millions of normalized zeros near the $10^{20}$\textsuperscript{th} and the $10^{22}$\textsuperscript{nd} zero of $\zeta(s)$. His computed correlation function shows remarkable agreement with Montgomery's form (see Figure \ref{fig:odlyzko}).

As we shall see, this work continues to have a profound impact on developments in contemporary Number Theory.

In the remaining two sections, we explore one statistic from random matrix theory (the density of eigenvalues) and one from number theory (the 1-level density of low-lying zeros). Though these statistics are not exactly analogous, they are similar. The reason we chose to study these two are that the general steps of the proofs are similar, and thus this provides a nice introduction to how intuitions and methods in one field can be transferred to another.


\section{Wigner's Semi-circle law}\label{sec:wigner}

\subsection{Wigner's Semi-circle Law (Statement)}\label{sec:wignsemicircstatement}

We state and prove a version of Wigner's semi-circle law below. We refer the reader to \cite{ERSY, ESY, TV1, TV2} for the most general version and proof of the semi-circle law as well as spacings between adjacent eigenvalues. We content ourselves with this special case as this version is easy to state and prove, and the conditions are frequently satisfied in practice.

\begin{thm}[Wigner's semi-circle law] Consider the ensemble of $N \times N$ real symmetric matrices with entries
independent, identically distributed random variables from a fixed probability distribution $p(x)$ with mean $0$, variance $1$, and other moments finite. Then
for almost all $A$, as $N \to \infty$
\begin{eqnarray}
\twocase{\mu_{A,N}(x)\  \longrightarrow \ }{ \frac{2}{\pi} \sqrt{1
- x^2}}{if $|x| \le 1$}{0}{otherwise.} \nonumber\
\end{eqnarray}
In other words, the number of normalized eigenvalues in an interval $[a,b] \subset [-1,1]$ is found by integrating the semi-circle over that interval. \end{thm}

Note that such a result could never hold for all $A$, as given any $\gep > 0$ there is always a small (though rapidly tending to zero!) probability that we've chosen a matrix that is within $\gep$ units from being a diagonal (i.e., each non-diagonal entry is at most $|\gep|$).

\begin{figure}
\begin{center}
\scalebox{.68}{\includegraphics{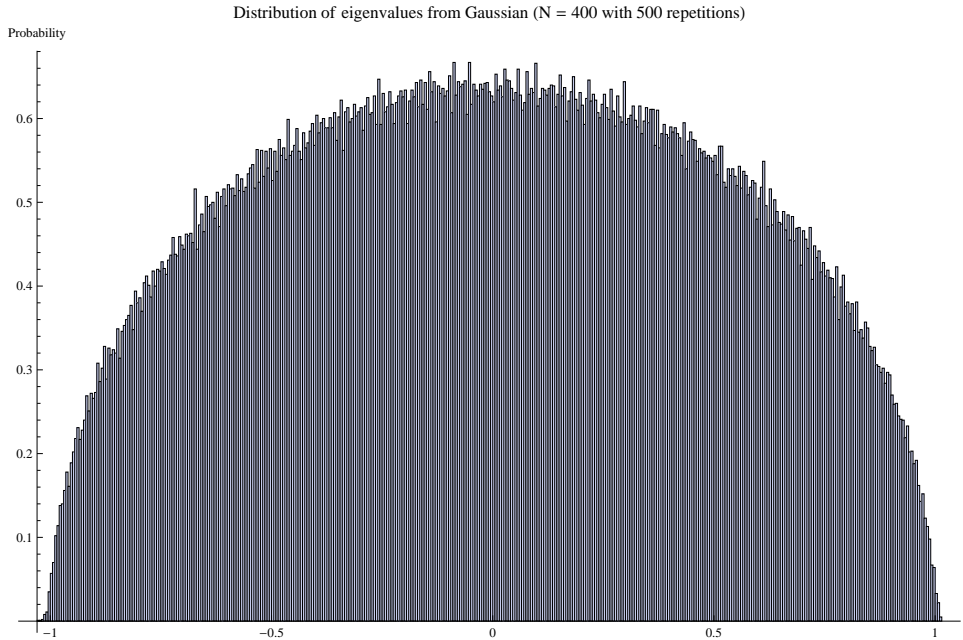}}
\caption{\label{fig:gaussian-400-500-a} A histogram plot of the normalized eigenvalues for 500 matrices, each $400 \times 400$. The entries are chosen independently from the standard normal $p(x) = (2\pi)^{-1/2} \exp(-x^2/2)$.}
\end{center}
\end{figure}

\begin{figure}
\begin{center}
\scalebox{.68}{\includegraphics{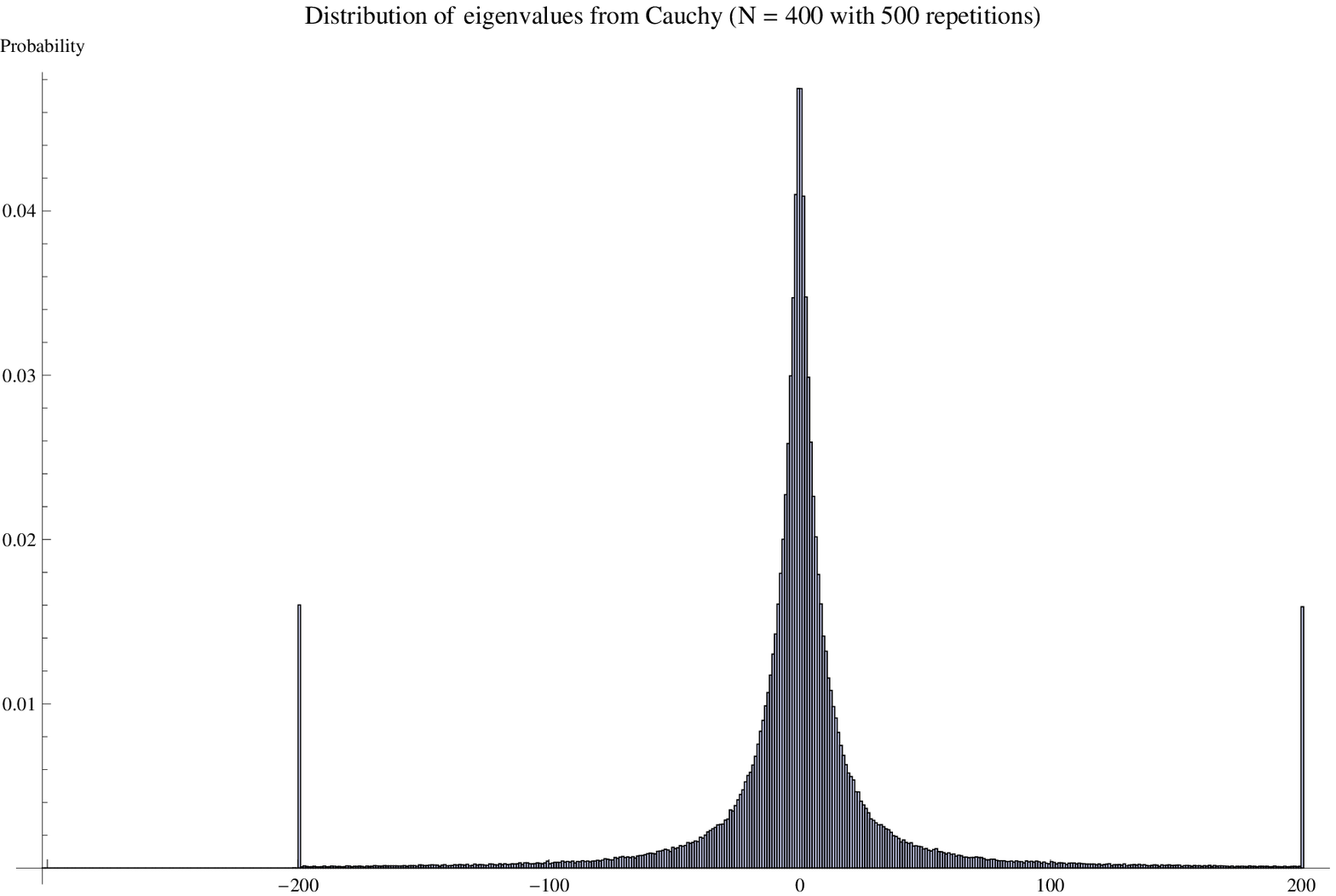}}
\caption{\label{fig:Dist2d} A histogram plot of the normalized eigenvalues for 500 matrices, each $400 \times 400$.  The entries are drawn from the Cauchy distribution $(\pi(1+x^2))^{-1}$. The bin on the extreme right represents all normalized eigenvalues that larger or large (and similarly for the bin on the extreme left).}
\end{center}
\end{figure}

For example, consider Figures \ref{fig:gaussian-400-500-a} and \ref{fig:Dist2d}. In the first we've drawn the entries from the standard normal. This satisfies the conditions of Wigner's semi-circle law, and we see already that with just $400 \times 400$ matrices the fit is excellent.

How essential are the conditions in the theorem? Does the result hold even if these conditions are violated, but perhaps the proof is just harder (or currently unknown)? To investigate this, we choose instead of the standard normal the Cauchy distribution $(\pi(1+x^2))^{-1}$. This distribution clearly has infinite variance, and thus obviously fails to satisfy the conditions. (It also has no mean as the integral of $|x|$ is infinite.) We see that the behavior is decidedly non-semi-circular (the huge probabilities at the end are the probabilities of observing an eigenvalue that far or further).

We will use the Method of Moments to prove the semi-circle law. We briefly summarize how we can pass from knowledge of the moments to knowledge of the eigenvalue distribution. Recall the $k$\textsuperscript{th} moment is $\sum_{i=1}^N \lambda_i(A)^k / 2^k N^{\frac{k}2+1}$. Imagine we had a $1\times 1$ matrix and we knew the first moment of the eigenvalues (well, here it would just be eigenvalue). We have one equation in one unknown: \be \lambda_1(A) / 2 \ = \ \mu_1; \ee this is clearly solvable and we can express $\lambda_1(A)$ in terms of $\mu_1$. Imagine now we have a $2 \times 2$ matrix. Then we have two equations in two unknowns: \bea \frac{1}{2^{5/2}} \left(\lambda_1(A) + \lambda_2(A)\right) & \ = \ & \mu_1 \nonumber\\ \frac{1}{2^4} \left(\lambda_1(A)^2 + \lambda_2(A)^2\right) & \ = \ & \mu_2. \eea For almost all $\mu_1$ and $\mu_2$ this is solvable (actually, we do not have to worry about this ever not being solvable, as the $\lambda_i(A)$ are always drawn from a matrix, and thus the equations will be consistent). We can therefore express the two eigenvalues in terms of the first two moments. Similarly, if we looked at the first three moments of a $3\times 3$ matrix we would have enough information to find the eigenvalues. In the general case, we need to know the first $N$ moments to find the eigenvalues of an $N\times N$ matrix;\footnote{Of course, one has to invert these relations to find the eigenvalues!} as we are letting $N\to\infty$, we need to compute all the moments to determine the eigenvalues.

The idea of the proof is as follows. For each matrix $A$ we calculate its moments; let us denote the $k$\textsuperscript{th} moment of $A$ by $M_{A,N}(k)$. The expected value of this is \be M_N(k) \ = \ \E[M_{A,N}(k)] \ = \ \int_{-\infty}^\infty \cdots \int_{-\infty}^\infty M_{A,N}(k) {\rm Prob}(A) dA. \ee We then show that $\lim_{N\to\infty} M_N(k) = C(k)$, the $k$\textsuperscript{th} moment of the semi-circle. This is almost, but not quite, enough to then conclude that a central limit theorem type situation occurs\footnote{See \cite{GS,Fel} for proofs of the central limit theorem, or \cite{MT-B} for a sketch of the proof.}, and a generic eigenvalue measure is close to the system average (which converges to the semi-circle as $N\to\infty$). The reason it is not sufficient is that we must also control the variances\footnote{For example, imagine for all $N$ we always had half the moments equal 0 and the other half equal $2C(k)$; then the average is $C(k)$ but no measure is close to the system average.}; however, this is easily done by similar arguments (see for example \cite{HM,MMS}).\\

\subsection{Wigner's Semi-circle Law (Sketch of Proof)}

We sketch the proof of Wigner's Semi-circle Law. As we've stated earlier, the reason we chose to prove this (as opposed to many of the other results) is that this proof is mostly self-contained, and highlights the key features of proofs in the subject.

There are typically three steps in working with random matrix ensembles (or the corresponding number theory quantities). We state these steps below, and then elaborate in great detail.

\ben

\item Determine the correct scale.\\

\item Develop an explicit formula relating what we want to study to something we understand. \\

\item Use an averaging formula to analyze the quantities above. \\

\een

Note it is not always trivial to figure out what is the correct statistic to study, and frequently very advanced combinatorics are needed to analyze the quantities.\footnote{Many of the papers in the field have large sections devoted to handling combinatorics; see for instance \cite{HM,Rub,RS}. Interestingly, sometimes the combinatorics \emph{cannot} be handled, and in \cite{Gao} we believe the number theory and random matrix theory agree where both have been calculated, but we cannot do the combinatorics to prove this.}

We describe these steps in detail for our random matrix ensembles and the semi-circle law. The key input is the following basic result from linear algebra, the Eigenvalue Trace Lemma.\footnote{The proof is trivial if $A$ is diagonal or upper diagonal, following by definition. As we only need this result for real symmetric matrices, which can be diagonalized, we only give the proof in this case. Let $S$ be such that $A = S\Lambda S^{-1}$ with $\Lambda$ diagonal. The claim follows from $A^k = S \Lambda^k S^{-1}$ and ${\rm Trace}(ABC) = {\rm Trace}(BCA)$.}

\begin{thm}[Eigenvalue Trace Lemma] Let $A$ be an $N\times N$ matrix with eigenvalues $\lambda_i(A)$. Then \bea {\rm Trace}(A^k) \ = \ \sum_{n=1}^N \lambda_i(A)^k. \eea As the trace of a matrix is the sum of its diagonal entries, \bea {\rm Trace}(A^k) \ = \ \sum_{i_1=1}^N \cdots \sum_{i_k = 1}^N a_{i_1i_2} a_{i_2i_3} \cdots a_{i_Ni_1}. \eea \end{thm}

The Eigenvalue Trace Lemma allows us to do the first two steps, namely determine the correct scale and relate what we want to study (the eigenvalues) to what we know (the matrix elements we randomly choose). We'll see later the analogue of this in number theory.\\

$\diamond$ For the first step, we take $k=2$ and find \be {\rm Trace}(A^2) \ = \ \sum_{i=1}^N \sum_{j=1}^N a_{ij} \ = \ \sum_{\ell=1}^N \lambda_\ell(A)^2 = N \langle \lambda(A)^2\rangle, \ee where $\langle \lambda(A)^2\rangle$ denotes the average of the square of the eigenvalues.  As $a_{ij} = a_{ji}$ and these are drawn from a probability distribution with mean 0 and variance 1, we have $\E[a_{ij}^2] = 1$, as this is just the variance. Thus the expected value of $N$ times the average eigenvalue square is just $N^2$, so the average eigenvalue square is of size $N$, so (heuristically) the average eigenvalue is of size $\sqrt{N}$.\footnote{With a little more work, we could calculate the variance of the average eigenvalue squared, using either the Central Limit Theorem or even Chebyshev's Theorem (from probability).} Why do we then normalize the eigenvalues by dividing by $2\sqrt{N}$ instead of $\sqrt{N}$? The reason is to make the final formula `clean' (i.e., this allows us to say the semi-circle law instead of the semi-ellipse); arguments such as these will capture the dependence on the key parameter (in this case, the $N$-dependence), but it will not catch constant dependence.\footnote{This is similar in some sense to dimensional analysis arguments in physics, which detect parameter dependence but not the constants. For example, imagine a pendulum of mass $m$ (in kilograms), length $L$ (in meters) where the difference in rest height (when the pendulum is down) and the raised height (where it is at angle $\theta$) is $L_0$ meters. We assume the only force acting is gravity, with constant $g$ (in meters per second squared). The period is how long (in seconds) it takes the pendulum to do a complete cycle, must be a function of $m$, $L$, $L_0$ and $g$; however, the only combinations of these quantities that give units of seconds are $\sqrt{L/g}$ and $\sqrt{L_0/g}$; thus the period must be a function of these two expressions. The correct answer turns out to be (at least for small initial displacements) approximately $2\pi \sqrt{L/g}$; we are able to deduce the correct functional form, though the constants are beyond such analysis.}\\

$\diamond$ For the second step, we want to understand the eigenvalues but it is the matrix elements we choose. The Eigenvalue Trace Lemma says $\sum \lambda_i(A)^k = {\rm Trace}(A^k)$; thus \eqref{eq:kthmomentmank} becomes \be\label{eq:kthmomentmankimproved} M_{A,N}(k) \ = \ \frac{\sum_{i=1}^N
\lambda_i(A)^k}{2^k N^{\frac{k}2+1}} \ = \ \frac{{\rm Trace}(A^k)}{2^k N^{\frac{k}2+1}}. \ee This is a terrific exchange. We have discussed how knowing the moments of the eigenvalues suffices to determine the eigenvalues; this allows us to express these moments in terms of the quantities we are choosing. \\

$\diamond$ For the third and final step, we note that in order for \eqref{eq:kthmomentmankimproved} to be useful we must be able to average it over all $A$ in our family; in other words, we must compute \bea M_N(k) & \ = \ &  \E[M_{A,N}(k)] \nonumber\\ & \ = \ & \int_{-\infty}^\infty \cdots \int_{-\infty}^\infty \frac{{\rm Trace}(A^k)}{2^k N^{\frac{k}2+1}} {\rm Prob}(A)dA \nonumber\\ & \ = \ & \int_{-\infty}^\infty \cdots \int_{-\infty}^\infty \frac{{\rm Trace}(A^k)}{2^k N^{\frac{k}2+1}} \prod_{1\le i \le j \le N} p(a_{ij}) da_{ij}. \eea The advantage is that ${\rm Trace}(A^k)$ is a polynomial in the matrix entries, and the integrals above can be readily evaluated.

For example, the integral for the average second moment is \bea \frac{1}{2^2 N^{2}}
\int_{-\infty}^\infty \cdots \int_{-\infty}^\infty \sum_{i=1}^N
\sum_{j=1}^N a_{ji}^2 \cdot p(a_{11})da_{11}\cdots p(a_{NN}) da_{NN}.
\nonumber \eea The integration factors as \bea \int_{a_{ij} =-\infty}^\infty a_{ij}^2
p(a_{ij}) da_{ij}\ \ \cdot \ \prod_{ { (k,l) \neq (i,j) \atop k <
l}} \int_{a_{kl} =-\infty}^\infty p(a_{kl}) da_{kl} \ = \ 1
\nonumber\ \eea (the first piece is 1 because it is a variance and thus is 1 by assumption, while the others are 1 as this is one of the defining properties of a probability distribution).

While the second moment calculation looks simple, the higher moment calculations require more involved computations and combinatorics\footnote{It is not hard to show the odd moments vanish by simple counting arguments. For the even moments, if the $a_{ij}$'s are not matched in pairs then there is negligible contribution as $N\to\infty$. The proof follows by counting how many tuples there can be with a given matching, and then comparing that to $N^{k/2+1}$ (the proofs are somewhat easier if our distribution is even). For example, as we are assuming our distribution has zero mean, if ever there was an $a_{i_\ell i_{\ell+1}}$ that was unmatched (so neither $(i_\ell,i_{\ell+1})$ or $(i_{\ell+1},i_\ell)$ occurs as the index in any other factor in the trace expansion, then the expectation of this term must vanish as each $a_{ij}$ is drawn from a mean zero distribution. The number of valid pairings of the $a_{ij}$'s (where everything is matched in pairs) is $(2k-1)!! = (2k-1)(2k-3)\cdots$, which is the $2k$\textsuperscript{th} moment of the standard normal; not ever matching contributes fully, though, and this is why the resulting moments are significantly smaller than the Gaussian's.}, in particular the Catalan numbers.\footnote{The Catalan numbers are $C_n = \frac1{n+1} \ncr{2n}{n}$. They arise in a variety of combinatorial problems; see for example \cite{Stan}.} The point is these computations can be done, and this analysis completes the proof (see \cite{MT-B} for the general arguments, and \cite{Leh} for the combinatorics calculation).\footnote{If the ensemble of matrices had a different symmetry, the start of the proof proceeds as above but the combinatorics changes. For example, looking at Real Symmetric Toeplitz matrices (matrices constant along diagonals) leads to very different combinatorics (and in fact a different density of states than the semi-circle); see \cite{HM,MMS}. If we looked at $d$-regular graphs, the combinatorics differs again, this time involving local trees; see \cite{McK}.}

\subsection{Additional statistics}

There are numerous other statistics we could investigate; the density of normalized eigenvalues is by no means the most natural, but it does highlight the general features. The more fundamental statistic is the spacings between adjacent normalized eigenvalues, and not their density (though the density is used to rescale to have mean spacing 1, allowing us to compare apples and apples). These spacings can either be attacked directly, or through the $n$-level correlations and combinatorics. It is conjectured that for our ensembles of real symmetric matrices, the spacings between normalized eigenvalues converges to a universal measure independent of $p$. This measure is approximately $(\pi/2)x \exp\left(-(\pi/4)x^2\right)$. Until very recently this was only known if the matrix elements were chosen from normal distributions; however, there has been great progress since the original version of this paper was written. L. Erd\"os, J. A. Ramirez, B. Schlein, T. Tao, V. Vu and H.-T. Yau \cite{ERSY,ESY,TV1,TV2} have removed this assumption and greatly generalized the class of matrices where the conjecture is known; the interested reader should see these papers for details. As the best result is constantly being improved, the interested reader should check the arxiv, \texttt{http://arxiv.org/}, for the current status.

\begin{figure}
\begin{center}
\scalebox{.45}{\includegraphics{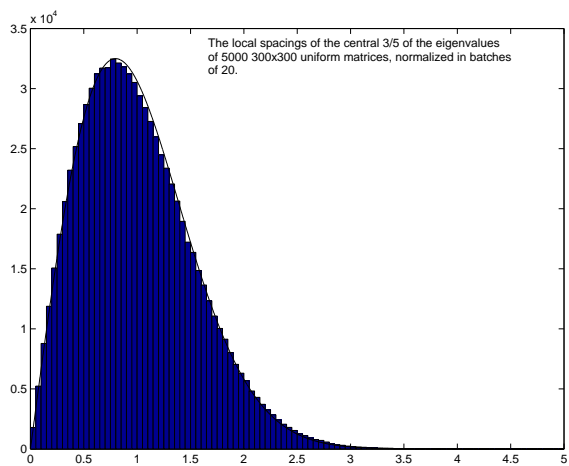}}
\caption{\label{fig:Dist1} Spacings between normalized eigenvalues of 5000 uniform matrices on $[-1,1]$ (size is 300).}
\end{center}
\end{figure}

Below we give some numerics from investigations in the bulk of the spectrum (i.e., normalized eigenvalues near 0).
Our first example is when $p$ is the uniform distribution on $[-1,1]$ (Figure \ref{fig:Dist1}). Already for $300\times 300$ matrices we see excellent agreement with the conjecture.

\begin{figure}
\scalebox{.3995}{\includegraphics{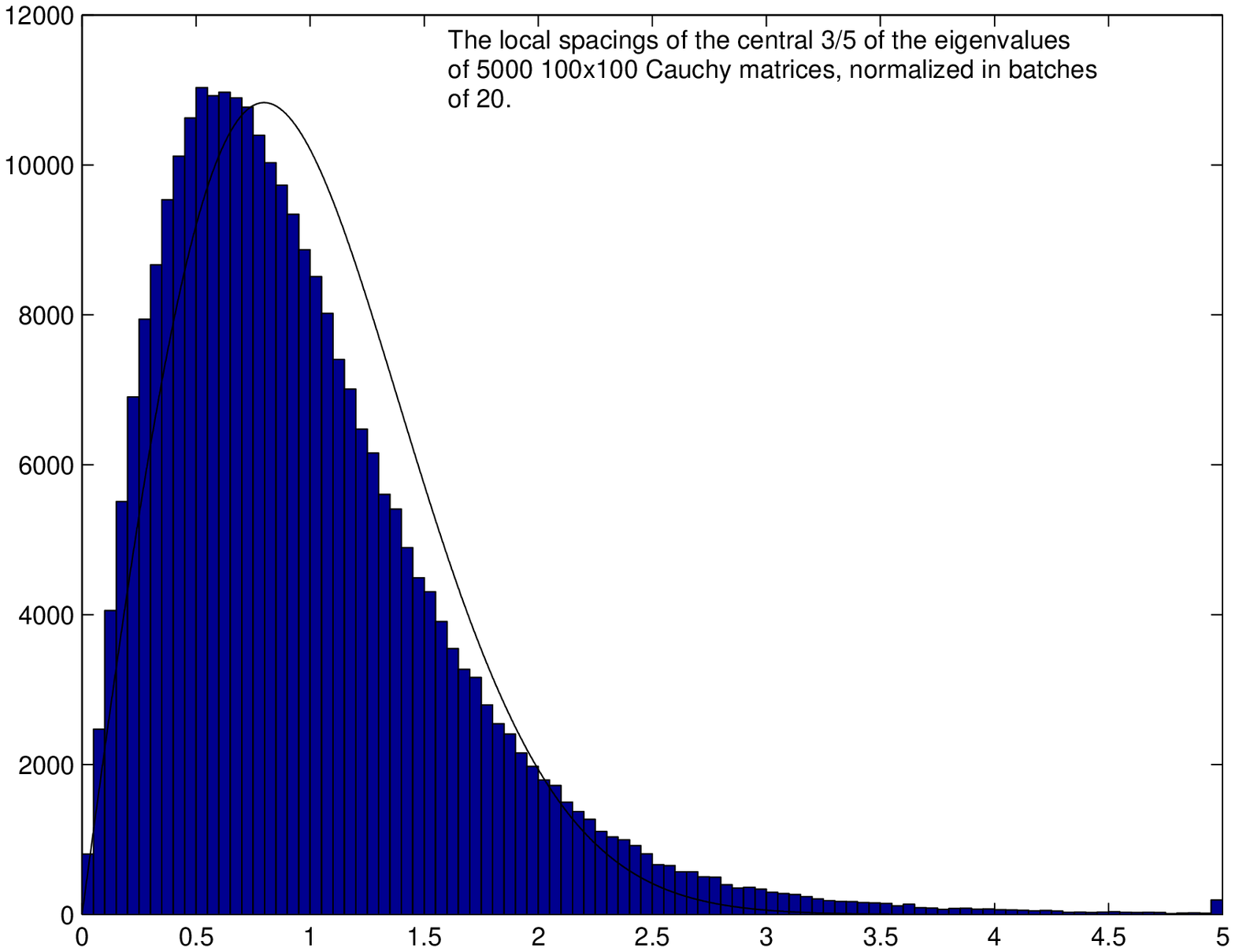}}\scalebox{.3995}{\includegraphics{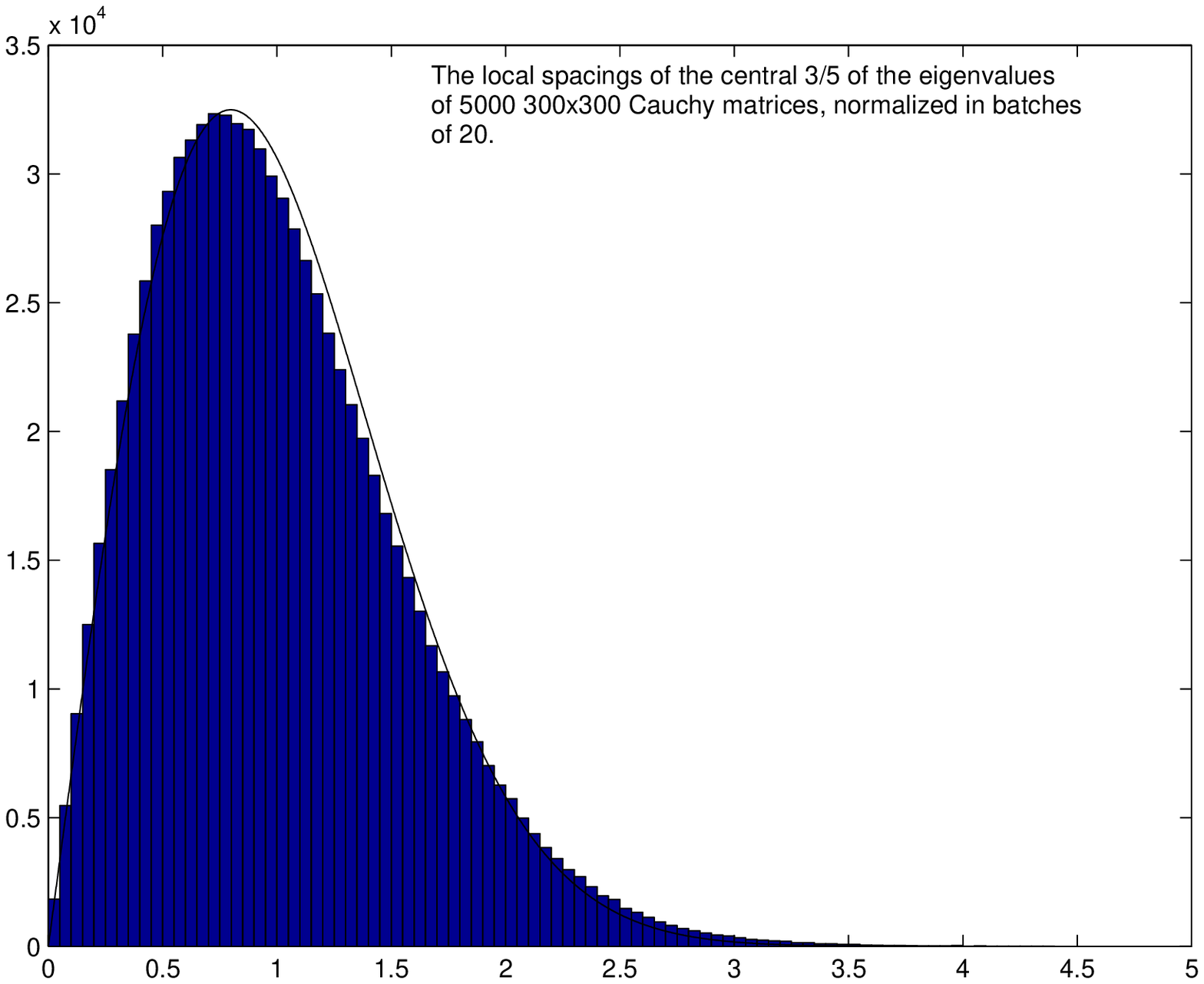}}  \\
\begin{center}
\caption{\label{fig:Dist2a} Spacings between normalized eigenvalues of 5000 Cauchy matrices; the left are $100\times 100$ and the right are $300\times 300$.}
\end{center}
\end{figure}

What about other distributions, for example the Cauchy density $(\pi(1+x^2))^{-1}$? We saw earlier that the density of states was decidedly non-semi-circular. It is a different story for the spacings (Figure \ref{fig:Dist2a}). Already for $300\times 300$ matrices we see excellent agreement with the conjecture.

There are numerous other ensembles where the density of states is non-semi-circular but the spacing between adjacent eigenvalues seems to agree with the conjecture. For example, McKay \cite{McK} proved the density of states of $d$-regular graphs\footnote{A graph is $d$-regular if each vertex is connected to exactly $d$ neighbors.} is Kesten's measure (which does converge to the semi-circle as $d\to\infty$), and simulations by many (including \cite{JMRR}) see the conjectured behavior. This is also apparent in more advanced tests, where the distribution of the largest eigenvalue is observed to follow a $\beta=1$ Tracy-Widom distribution (see \cite{MNS}). These distributions govern the largest eigenvalues in many settings; see \cite{TW1,TW2,TW3}. If the ensemble has a very different structure, however (such as the Toeplitz ensembles in \cite{BDJ,HM,MMS}) then both statistics could behave differently.


\section{From Random Matrix Theory to Number Theory}\label{sec:rmttonumb}

We now discuss the path from random matrix theory to number theory. As this story has been told numerous times, we concentrate on some illuminating aspects and refer the reader to the references previously mentioned. We concentrate on a very small subset of statistics and connections; for example, we will almost completely ignore the contributions to studying moments of the zeta function. Our goal is to explain how the behavior of some key statistics in number theory are the same as the corresponding statistics in random matrix theory. We concentrate on $\zeta(s)$ and its simplest generalization, Dirichlet $L$-functions, though these results hold for a larger class of $L$-functions as well (exactly how large a class is the subject of many research programs).

The starting point of our analysis is Riemann's Explicit Formula, which is a natural generalization of \eqref{eq:logderiv}. Let $\phi(s)$ be a `nice' function. We have \be \frac1{2\pi i} \int_{(c)} -\frac{\zeta'(s)}{\zeta(s)}\ \phi(s) ds \ = \  \frac1{2\pi i}\int_{(c)} \sum_{n=1}^\infty \Lambda(n) \frac{\phi(s)}{n^s}\ ds. \ee Shifting contours, by the residue theorem the left hand size is basically $\phi(1) - \sum_\rho \phi(\rho)$, while the right hand side is (for $s = \sigma + it$): \bea \int_{(c)} \sum_{n=1}^\infty \Lambda(n) \frac{\phi(s)}{n^s}\ ds & \ = \ & \sum_{n=1}^\infty \frac{\Lambda(n)}{n^\sigma} \int_{-\infty}^\infty \phi\left(\sigma + it\right) e^{-it \log n} dt.\ \ \eea The integral is basically the Fourier transform\footnote{The Fourier transform of $g$ is $\widehat{g}(\xi) = \int_{-\infty}^\infty g(x) e^{-2\pi i x \xi}dx$.} (up to some constants) of $f(t) = \phi(\sigma + it)$. A careful analysis \cite{MT-B} gives the following explicit formula relating sums of a test function over zeros of $\zeta(s)$ to sums of the Fourier transform over primes.

\begin{thm}[Explicit Formula]
Let $\sum_\rho$ denote the sum over the non-trivial zeros of $\zeta(s)$ (i.e., the zeros in
critical strip), $g$ an even Schwartz function\footnote{This means for any $m, n \ge 0$ that $\lim_{|x|\to\infty} (1+x^2)^m g^{(n)}(x)=0$ (i.e., $g$ and all its derivatives tend to zero faster than any polynomial).} of compact
support and $\phi(r) = \int_{-\infty}^\infty g(u) e^{iru}du$. Write $\rho$ as $1/2 + i \gamma_\rho$; if the Riemann Hypothesis is true then $\gamma_\rho$ is real. We have
\begin{eqnarray} & & \sum_\rho \phi(\gamma_\rho)   \ = \ 2 \phi\lp\frac{i}{2}\rp
-  \sum_{p} \sum_{k=1}^{\infty} \frac{2 \log p}{p^{k/2}}
g \left( k \log p \right)\nonumber\\
&& \ \ \ \ \ \ \ + \ \frac{1}{\pi} \int_{-\infty}^\infty \left( \frac{1}{iy -
\foh} + \frac{\Gamma'(\frac{iy}{2} +
\frac{5}{4})}{\Gamma(\frac{iy}{2} + \frac{5}{4})} - \foh \log \pi
\right) \phi(y) \, dy; \nonumber\
\end{eqnarray} note up to scale that $g$ and $\phi$ are essentially a Fourier transform pair.
\end{thm}

\subsection{Preliminaries}

\begin{figure}
\begin{center}
\scalebox{1}{\includegraphics{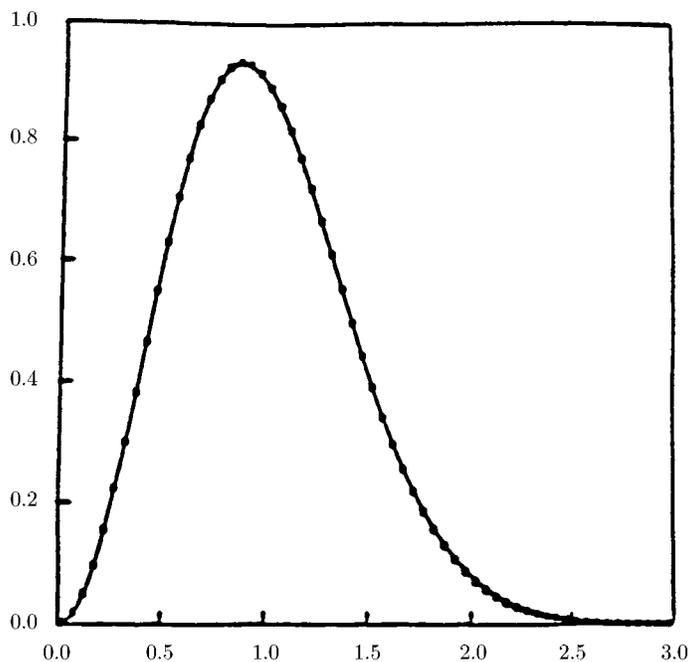}}
\caption{\label{fig:odlyzko2} 70 million spacings between adjacent zeros of $\zeta(s)$,
starting at the $10^{20}$\textsuperscript{th} zero, versus the corresponding results for eigenvalues of complex Hermitian matrices (from Odlyzko).}
\end{center}
\end{figure}

We now come to the key moment of our story, when Montgomery and Dyson \cite{Mon} noticed the agreement between the pair correlation of zeros of $\zeta(s)$ and eigenvalues of complex Hermitian matrices. The pair correlation statistic of a set $\{x_1, x_2, \dots\}$ is \be \lim_{N\to\infty} \frac{\#\{i \neq j \le N: x_i - x_j \in I\}}{N} \ee where $I$ is an arbitrary interval. We can generalize this to triple correlation (which would be how often pairs of differences are in a box) and higher; knowing all the correlations is equivalent to knowing the spacing between adjacent elements. Instead of using a box or hypercube we can use a smooth test function.\footnote{We want (1) $f(x_1, \dots, x_n)$ is symmetric; (2)
$f(x + t(1, \dots, 1)) = f(x)$ for $t\in \R$; (3) $f(x) \to 0$ rapidly as $|x| \to \infty$ in the hyperplane $\sum_k x_j =0$; see \cite{RS}.} Odlyzko \cite{Od1,Od2} observed phenomenal agreement between adjacent zeros and the corresponding distribution for spacings between adjacent eigenvalues of complex Hermitian matrices; see Figure \ref{fig:odlyzko2}.

Hejhal \cite{Hej} proved (for suitable test functions) that the triple correlation of $\zeta(s)$ agrees with random matrix theory, while Rudnick and Sarnak \cite{RS} showed agreement with the $n$-level correlations of any $L$-function attached to a cuspidal automorphic representation of ${\rm GL}_m/\Q$. To describe these $L$-functions in detail would be too much of a digression, so we will content ourselves with a very brief introduction, referring the interested reader to \cite{RS} for details. The Riemann zeta function \be \zeta(s) \ = \ \sum_{n=1}^\infty \frac{1}{n^s} \ = \ \prod_{p\ {\rm prime}} \left(1 - \frac1{p^s}\right)^{-1} \ee has an Euler product, a functional equation, and is conjectured to have all of its zeros in the critical strip $0 \le \Re(s) \le 1$ on the line $\Re(s) = 1/2$. The generalization is a series such as $\sum_{n=1}^\infty a_n / n^s$, where the $a_n$ are of arithmetic interest. In order to call this series an $L$-function, we require it to have certain properties (such as an Euler product and a functional equation, as well as some growth rates on the $a_n$'s). We call these $L$-functions, and they arise throughout number theory.\footnote{The earliest occurrence was probably in Dirichlet's work, who used $L$-functions attached to characters on $(\Z/m\Z)^\ast$ to study primes in arithmetic progressions. Another example are elliptic curve $L$-functions, which (at least conjecturally) give information about the group of Mordell-Weil group of rational solutions of the elliptic curve.}

In random matrix theory, in order to understand the behavior of the eigenvalues of one matrix $A$ we embedded it in a family of random matrices, and showed that with high probability the behavior of the eigenvalues of $A$ are close to the ensemble average (at least as $N\to\infty$). For the $n$-level correlations, we do not need to perform any averaging over $L$-functions; we may study an individual $L$-function. The reason is that the density of zeros in the critical strip whose imaginary part is of size $T$ is (up to some constants) of size $1/\log T$; in other words, the higher up we go, the more densely the zeros of an $L$-function are packed together, and thus one $L$-function provides enough zeros high up to average.

The results mentioned above suggested that, for the purposes of number theory, it sufficed to know how random complex Hermitian matrices behave, as the zeros of $\zeta(s)$ (and other $L$-functions) high up on the critical line showed remarkable agreement with these eigenvalues. This turned out, however, to only be part of the story. The reason is that the $n$-level correlations are insensitive to finitely many zeros. In other words, if we were to remove the 1701 zeros nearest to the critical point $s=1/2$, the $n$-level correlations would not change.\footnote{This is because the zeros are tending to infinity. Thus, given any zero and any finite box, only finitely many zeros can be associated to it such that the required differences lie in the box. Therefore, this zero has negligible contribution in the limit as we are dividing by $N$.} This is a major problem for number theory, as often we expect there to be behavior at the central point of arithmetic interest.\footnote{For example, the Birch and Swinnerton-Dyer conjecture states the order of vanishing at the central point of an elliptic curve $L$-function equals the rank of the Mordell-Weil group of rational solutions of the elliptic curve; this is quite important information which we do not wish to discard!}

Katz and Sarnak \cite{KaSa1,KaSa2} showed that, as the size of the matrices tends to infinity, the $n$-level correlations of complex Hermitian matrices also equals those of $N\times N$ unitary matrices, as well as its orthogonal and symplectic subgroups.\footnote{These classical compact groups are much more natural random matrix ensembles. In our original formulation, we chose the matrix elements randomly and independently from some probability distribution $p$. What should we take for $p$? The GOE and GUE ensembles, where the entries are chosen from Gaussians, arise by imposing invariance on ${\rm Prob}(A)dA$ under orthogonal (respectively unitary) transformations; these are natural conditions to impose as the probability of a transformation should be independent of the coordinate system used to write it down. The classical compact groups come endowed with a natural probability, namely Haar measure.} Thus when we say that the zeros behave like eigenvalues of complex Hermitian matrices, we could have also said they behave like eigenvalues of unitary matrices (or one of their subgroups).\footnote{The eigenvalues of a unitary matrix are of the form $e^{i\theta}$. To see this, let $v$ be an eigenvector of $U$ with eigenvalue $\lambda$. Note $v^\ast U^\ast U v = v^\ast v$, which gives $|\lambda|^2 ||v||^2 = ||v||^2$, so $|\lambda| = 1$. Thus, similar to real symmetric and complex Hermitian matrices, we can parametrize the eigenvalues by a real quantity.} We thus need a new statistic to study which will `break' this symmetry, and say which ensemble is truly modeling the behavior. Further, our statistic should take into account the behavior near the central point, as interesting arithmetic occurs there. One popular choice is the $1$-level density, which we now describe.

\subsection{$1$-level Density (Preliminaries)}

Let $\phi(x)$ be an even Schwartz function. This means for any $m, n \ge 0$ that $\lim_{|x|\to\infty} (1+x^2)^m \phi^{(n)}(x)$ (i.e., $\phi$ and all its derivatives tend to zero faster than any polynomial). We also assume the Fourier Transform\footnote{We use the normalization $\widehat{\phi}(\xi) = \int_{-\infty}^\infty \phi(x) e^{-2\pi i x \xi}dx$. The Fourier transform has many nice properties on the Schwartz space (see \cite{SS} for example).} is compactly supported; this means there is some $\sigma < \infty$ such that $\widehat{\phi}(\xi) = 0$ if $|\xi| \ge \sigma$.

Consider an $L$-function \be L(s,f) \ = \ \sum_{n=1}^\infty \frac{\lambda_f(n)}{n^s}; \ee we assume this series converges for $\Re(s) > 1$, has a meromorphic extension to all of $\C$ satisfying a functional equation, and has an Euler product.\footnote{These are very strong conditions, and most choices of $\lambda_f(n)$ will not satisfy these requirements. Fortunately there are many choices that do work, and these frequently encode information about arithmetically interesting problems. The most studied examples include Dirichlet $L$-functions, modular and Maass forms; see \cite{IK} for details.} We have remarked that high up, the spacing between zeros is like $1/\log T$ at height $T$; what is it near $s=1/2$? The answer can be deduced from an analysis of the functional equation, which shows there is some number $C_f$ (called the analytic conductor) such that the zeros near the central point are spaced on the order of $1/\log C_f$. This suggests we study the following statistic:

\begin{defi}[The $1$-level density] Let $\phi$ be an even Schwartz function and $L(s,f)$ an $L$-function as above. The $1$-level density is
\begin{eqnarray}
D_{1, f}(\phi) \ &=& \ \sum_{j}
\phi\left(\gamma_{f,j} \frac{\log C_f}{2\pi}\right).
\end{eqnarray}
\end{defi}

We may generalize the above to $n$-level densities; while for some applications it is essential to understand these generalizations\footnote{For example, there are three flavors of orthogonal symmetry, and their $1$-level densities are indistinguishable if the support of $\widehat{\phi}$ is contained in $(-1,1)$; however, the $2$-level densities of all three are distinguishable for arbitrarily small support \cite{Mil1}. Another example is in obtaining better decay rates for high vanishing at the central point in a family \cite{HM}.}, for many purposes studying the $1$-level density suffices. This statistic differs in several important ways from the $n$-level correlation.

The first is that individual zeros now contribute in the limit. Moreover, most of the contribution is from the zeros near the central point (thus this statistic is sensitive to what is happening there). This is because $\phi$ is of rapid decay, so once we are a couple of average spacings away, there is negligible contribution. There is a trade-off, namely it no longer suffices to study just one $L$-function. The reason is we always need something to average over, and there are just too few zeros near the central point on this scale. The solution is to look at the zeros near the central point for many $L$-functions that share common properties. We average the $1$-level densities over the family. Unlike the $n$-level correlations, where we are looking high up on the critical line and see the same behavior in all $L$-functions, we see very different behavior near the central point, depending on what family of $L$-functions we study.

Katz and Sarnak \cite{KaSa1,KaSa2} conjecture that to any `nice' family of $L$-functions, as the conductors tend to infinity the $1$-level density of the family agrees with the $N\to\infty$ scaling limit of a classical compact group (typically $N\times N$ unitary, orthogonal or symplectic matrices). Moreover, these groups \emph{all} have distinguishable behavior. In other words, the universality seen in the $n$-level correlations is broken.

Before describing the proof, we give some examples of families of $L$-functions and the corresponding symmetries.

\ben

\item Dirichlet $L$-functions: Let $m$ be a prime and consider all non-principal\footnote{This means we avoid the trivial character $\chi_0(n) = 1$ if $n$ is relatively prime to $m$ and $0$ otherwise; this character gives rise to a simple modification of $\zeta(s)$.} Dirichlet characters $\chi$ from $(\Z/m\Z)^\ast$ to the complex numbers of absolute value 1. To each character $\chi$ we have an $L$-function $L(s,\chi) = \sum_n \chi(n)/n^s = \prod_p(1 - \chi(p)p^{-s})^{-1}$. As $q\to\infty$, the behavior agrees with the scaling limit of unitary matrices.\footnote{We could consider the related family of quadratic Dirichlet characters coming from a fundamental discriminant $d \in [X,2X]$ with $X\to\infty$; this family agrees with the scaling limit of symplectic matrices.}\\

\item Cuspidal newforms: Let
\be \Gamma_0(N) \ = \ \left\{ \mattwo{a}{b}{c}{d}: { ad-bc \ = \ 1
\atop c \equiv 0 (N) } \right\}. \ee We say $f$ is a weight $k$ holomorphic cuspform of level $N$ if \be \forall \gamma \in \Gamma_0(N), \ \ f(\gamma z) \ = \ (cz+d)^k
f(z), \ee where $\gamma z = (az+b)/(cz+d)$. As $k$ or $N$ tend to infinity, the behavior agrees with the scaling limit of orthogonal matrices.\footnote{There are actually three flavors of orthogonal groups. If all the signs of the functional equation are even, the corresponding group should be ${\rm SO}(2N)$, and similarly if the signs are all odd. We refer the reader to the previously mentioned surveys for the details.} \\

\item Elliptic curves: Let $\mathcal{E}: y^2 = x^3 + A(T)x + B(T)$ be an elliptic curve over $\Q(T)$. For each $t\in \Z$ we can specialize and get an elliptic curve $E_t: y^2 = x^3 + A(t)x+B(t)$. We can build an $L$-function, where $a_p(t)$ is related to the number of solutions to $y^2 \equiv x^3 + A(t)x+B(t) \bmod p$. Our family is now $t \in [X,2X]$ with $X\to\infty$, and these families have orthogonal symmetry.    \\

\een

In these and many other cases (see \cite{DM1,FI,Gao,Gu,HM,HR,ILS,KaSa2,Mil1,Mil5,Ro,Rub,Yo2} for a representative sampling of results), we can show for suitably restricted $\phi$ that the $1$-level density agrees with the scaling limit of one of the mentioned classical compact groups.

\subsection{$1$-level Density (Proofs)}

We briefly describe how the proofs proceed. We concentrate on the family of Dirichlet characters with prime conductor $q$ tending to infinity. As in the proof of Wigner's semi-circle law, there are three steps.\\

$\diamond$ We first must determine the correct scale to study the zeros near the central point. The answer can be shown to follow from the functional equation; in this case, we normalize the zeros by the factor $(\log (m/\pi))/ (2\pi)$. \\

$\diamond$ The next step is to relate what we want to study, namely the zeros near the central point, to something we have a chance of understanding, in this case the coefficients of the $L$-function, $\chi(n)$. We do this with a generalization of Riemann's explicit formula (this is the analogue of the Eigenvalue Trace Lemma). Our result is a straightforward generalization of the formula Riemann used to connect the zeros of $\zeta(s)$ with the distribution of primes. The difference here is that we have a more general test function. It can be shown (see \cite{ILS,RS}) that
for $\phi$ an even Schwartz
function and $L(s,\chi) = \sum_n \chi(n)/n^s$ a Dirichlet $L$-function from a
non-trivial character $\chi$ with conductor $m$ and zeros $\rho =
\foh + i\gamma_{\chi,\rho}$, then \bea & & \sum_\rho
\phi\left(\gamma_{\chi,\rho} \frac{ \log(m/\pi)}{2\pi}\right) \ = \ \int_{-\infty}^{\infty}
\phi(y)dy \nonumber\\ & & -2 \sum_p \logpm
\widehat{\phi}\left(\logpm\right) \frac{\chi(p)}{p^{1/2}} \nonumber\\ & & -2 \sum_p \logpm \widehat{\phi}\left(\logptm\right) \frac{\chi^2(p)}{p} + O\left(\frac1{\log m}\right).  \eea Note the left hand side is a sum over zeros and the right hand side a sum over the coefficients in the $L$-function. We also have $\widehat{\phi}$ on the right hand side. We now see how the support condition enters. As we assume $\widehat{\phi}$ is supported in $(-\sigma, \sigma)$, this restricts the sums on the right to having only finitely many terms. This is the $1$-level density for one Dirichlet $L$-function $L(s,\chi)$. We now average over all non-principal $\chi$ (i.e., $\chi \neq \chi_0$). We will see below that there are $m-2$ such characters, and thus we obtain \bea & & \frac1{m-2} \sum_{\chi \neq \chi_0} \sum_\rho
\phi\left(\gamma_{\chi,\rho} \frac{ \log(m/\pi)}{2\pi}\right) \ = \ \int_{-\infty}^{\infty}
\phi(y)dy \nonumber\\ & & -2 \sum_p \frac{\log p}{p^{1/2}\log(m/\pi)}
\widehat{\phi}\left( \logpm\right) \frac1{m-2}\sum_{\chi \neq \chi_0} \chi(p) \nonumber\\ & & -2 \sum_p \frac{\log p}{p\log(m/\pi)} \widehat{\phi}\left(\logptm\right)\frac1{m-2} \sum_{\chi \neq \chi_0}\chi^2(p) + O\left(\frac1{\log m}\right). \nonumber\\ \eea \\

$\diamond$ Similar to the proof of Wigner's semi-circle law, our explicit formula would be useless unless we can perform the averaging over the family. We briefly review some needed results about these characters (see for instance \cite{MT-B}), and then show how to handle the sums. Unfortunately, the averaging formulas in number theory are significantly worse than the corresponding averaging formulas in random matrix theory; this leads to far more restricted results in number theory.\footnote{In fact, the case being studied here has the best averaging formula! For cuspidal newforms we have the Petersson formula, and for families of elliptic curves we can use periodicity in evaluating Legendre sums of cubic polynomials. In general, however, unless our family is obtained in some manner from a family such as one of these, we do not possess the needed averaging formula.}

For $m$ prime, the group $(\Z / m\Z)^*$ is cyclic of order $m-1$;\footnote{$(\Z/m\Z)^\ast$ is the set $\{1,2,\dots,m-1\}$ where multiplication and addition are modulo $m$; for example, if $m=17$, $x=11$ and $y=9$ then $xy = 99 = 5 \cdot 17 + 14$, so $xy \equiv 14 \bmod 17$, while $x+y = 20 \equiv 3 \bmod 17$. The hardest step in proving that this is a group under multiplication is finding inverses; one way to accomplish this is with the Euclidean algorithm.} denote its generator by $g$. Let $\zeta_{m-1} = e^{2 \pi i /(m-1)}$. The principal character $\chi_0$ is given by
\begin{equation}
\twocase{\chi_0(k)\ =\ }{1}{if $(k,m) = 1$}{0}{if $(k,m) > 1$.}
\end{equation} As the characters are group homomorphisms from $(\Z/m\Z)^\ast$ to complex numbers of absolute value 1, the $m-2$ primitive characters are determined by their action on $g$.\footnote{These properties mean $\chi(1) = 1$, $\chi(xy) = \chi(x)\chi(y)$ and $\chi(x^r) = \chi(x)^r$. Further, though initially only defined on $(\Z/m\Z)^\ast$, we extend the definition to all of $\Z$ by setting $\chi(n+\ell m) = \chi(n)$. As $g^{m-1} \equiv 1 \bmod m$, $\chi(g^{m-1}) = 1$ or $\chi(g)^{m-1} = 1$ for all $\chi$. Thus $\chi(g)$ has to be an $m-1$ root of unity. We see each of these roots gives rise to a character (one of which is the principal or trivial character); by multiplicativity once we know the character's action on the generator we know it everywhere.}
Thus there exists an $\ell$ such that $\chi(g) = \zeta_{m-1}^\ell.$ A simple calculation (use the explicit representation for the characters and the geometric series formula) shows that
\begin{equation}
\twocase{\sum_\chi \chi(k)\ =\ }{m-1}{if $k \equiv 1 \bmod m$}{0}{otherwise,}
\end{equation} where we are summing over \emph{all} characters, including the principal one. It is easy to remove the contribution from the principal character, and we find
for any prime $p \neq m$
\be \twocase{\sum_{\chi \neq \chi_0} \chi(p) \ = \ }{-1+m-1}{$p
\equiv 1 (m)$}{-1}{otherwise.}  \ee
This is the desired averaging formula. We substitute it into
\be  \sum_p \frac{\log p}{p^{1/2}\log (m/\pi)}\
\hphi\left(\logpm\right) \fommt \sum_{\chi \neq \chi_0}\chi(p).\ee We write $f(x) \ll g(x)$ if there is a $C$ such that for all $x$ sufficiently large, $|f(x)| \le C g(x)$. Our function $\phi$ is bounded, and as $p \le m^\sigma$, $\log p \ll \sigma \log (m/\pi)$ as $m\to\infty$.
A simple calculation shows there is no contribution if $\sigma < 2$: \begin{eqnarray}
& & \frac{-2}{m-2} \sum_p^{m^\sigma} \logpm \hphi\left(\logpm\right)
\cdot p^{-\foh} \nonumber\\  & + & 2\frac{m-1}{m-2} \sum_{p \equiv 1 (m)}
^{m^\sigma} \logpm \hphi\left(\logpm\right)\cdot p^{-\foh} \nonumber\\ &
\ll & \frac{1}{m}\sum_p^{m^\sigma} p^{-\foh} \ + \ \sum_{p \equiv
1 (m)}^{m^\sigma} p^{-\foh} \ll
\frac{1}{m}\sum_k^{m^\sigma} k^{-\foh} \ + \ \sum_{ {k \equiv 1
(m) \atop k \ge m+1} }^{m^\sigma} k^{-\foh}  \nonumber\\ & \ll &
\frac{1}{m}\sum_k^{m^\sigma} k^{-\foh} \ + \
\frac{1}{m}\sum_{k}^{m^\sigma} k^{-\foh} \ll
\frac{1}{m} m^{\sigma/2}.
\end{eqnarray}
It is conjectured that there should be no contribution for any finite $\sigma$;\footnote{This is not surprising, as the above argument is quite crude, where we have inserted absolute values and thus lost the expected cancelation due to the terms having opposite sign.} extending this further is related to some of the deepest questions about how primes are distributed in arithmetic progressions.

\subsection{Nuclear Physics Interpretation}

It is interesting to interpret our results on the $1$-level density in the language of nuclear physics:\\

\bea  \mbox{Zeros of $L$-functions} &  \ \longleftrightarrow \ & \mbox{Energy levels of heavy nuclei} \nonumber\\
& & \nonumber\\ \mbox{Schwartz\ test\ function} & \longleftrightarrow & \mbox{Neutron} \nonumber\\
& & \nonumber\\ \mbox{Support\ of\ test\ function} & \longleftrightarrow &
\mbox{Neutron Energy}. \nonumber\ \eea \ \\

We expand on the above. Assuming the (Generalized) Riemann Hypothesis\footnote{This just asserts that any `nice' $L$-function has all of its zeros in the critical strip having real part equal to 1/2.}, the zeros of our $L$-functions lie on the critical line $\Re(s) = 1/2$. Thus we may order them, and it makes sense to talk about spacings between adjacent zeros. Following the P\'{o}lya-Hilbert dream (and there are numerous people pursuing this), we may even try to search for a physical system whose energy levels are these zeros! Regardless, we have two sequences of real numbers: the zeros of our $L$-function(s) and the energy levels of our heavy nucleus (nuclei).

How do we understand the structure of the nucleus? We bombard it with low energy neutrons, and see what happens. The analogy on the number theory side is we `bombard' the zeros of our $L$-function with our Schwartz test function; what we `see' now is a sum over primes.

In physics, ideally we would like to be able to send in a neutron with any energy; unfortunately, current technology only allows us to send in neutrons with energy in a given band. Thus we cannot obtain perfect information about the internal structure of the nucleus and its energy levels. This corresponds exactly with number theory, where the support of the Fourier transform of our Schwartz test function is playing the role of the neutron's energy. Bombarding the zeros with a test function is equivalent to summing the Fourier transform against related quantities, and our averaging formulas are only able to handle certain restricted sums.

It is worth dwelling on this last observation a little more. The Heisenberg Uncertainty Principle can be recast in mathematical terms as a statement about a function and its Fourier transform, namely it is not possible to simultaneously localize $f$ and $\widehat{f}$ (i.e., the product of the variances, their spreads about their means, cannot be too small). Ideally we would like to take $\delta(x-a)$ as our test function, as this would allow us to understand whether or not there are zeros at $a$. In particular, if we take $a=0$ we would understand the behavior at the critical point. Unfortunately the Fourier transform of $\delta(x)$ is identically 1; this corresponds to having absolutely \emph{no} control on the prime sum side.

\subsection{Future avenues}

Random matrix theory has enjoyed remarkable success in suggesting questions and predicting answers for number theory. The $n$-level correlations and densities are two of many examples. One problem, however, is that random matrix theory often cannot detect the arithmetic of the $L$-functions, and this must be added (in a sometimes unsatisfying manner). A terrific example of this is in the study of moments of $L$-functions (see \cite{CGo,CGh,CFKRS,KeSn1,KeSn2} and the references therein); see also the discussion below on the hybrid product formulas of Gonek, Hughes and Keating.

For families of $L$-functions $\{L(s,f_i)\}_{f_i\in\mathcal{F}_{i}}$ ($i \in \{1,2,\dots,I\}$), we can form the Rankin-Selberg convolution and study the family \begin{equation} \{L(s,f_1 \otimes \cdots \otimes f_I)\}_{(f_1,\dots,f_I) \in \mathcal{F}_1 \times \cdots \times \mathcal{F}_I}; \nonumber\ \end{equation} see \cite{IK} for details.\footnote{If $L(s,f_i)$ $=$ $L_{\infty;i}(s)$ $\prod_p$ $\prod_{j=1}^{m_i}$ $(1 - \alpha_{j;i}(p)p^{-s})^{-1}$, then $L(s,f_1 \otimes f_2)$ is related to a product over primes of $\prod_{j=1}^{m_1} \prod_{k=1}^{m_2}$ $(1 - \alpha_{j;1}(p)\alpha_{k;2}(p) p^{-s})^{-1}$. It is conjectured that these functions have functional equations, satisfy the Riemann hypothesis, et cetera. The existence of the Rankin-Selberg convolution is known for just a few choices of the $f_i$'s.} Given a family $\{L(s,f_i)\}_{f_i\in\mathcal{F}_{i}}$, the Katz-Sarnak conjecture states the behavior of zeros near the central point (as the conductors tend to infinity) agrees with the $N\to\infty$ scaling limit of a subgroup of unitary matrices $U(N)$. A natural question to ask is how the behavior of zeros in the family of convolutions is related to the behavior of the constituent families. Due\~nez-Miller show that if the family of $L$-functions are `nice' (see \cite{DM2} for statements and proofs\footnote{There are several difficulties with the proofs in general, ranging from not knowing properties of the Rankin-Selberg convolution in general to not having a good averaging formula over general families.}), then we can attach a symmetry constant $c(\mathcal{F}_{i})$ to each family satisfying the following conditions:
\ben

\item $c(\mathcal{F}_i)$ is 0 if the family has unitary symmetry, $1$ if the family has symplectic symmetry and $-1$ if the family has orthogonal symmetry;

\item $c(\mathcal{F}_i \times \mathcal{F}_j) = c(\mathcal{F}_i)\times c(\mathcal{F}_j)$.

\een In other words, for many families the symmetry type of the convolution is the product of the symmetry types.\footnote{A special case of this theorem was discovered by Due\~nez-Miller in \cite{DM1} in studies of a family of ${\rm GL}(4)$ and a family of ${\rm GL}(6)$ $L$-functions. The analysis there led to a disproof of a folklore conjecture that the theory of low-lying zeros of $L$-functions is equivalent to a theory of the distribution of signs of the functional equations in the family; see \cite{DM1,DM2} for details. The key ingredient in the proofs is the universality of the second moments of the Satake parameters $\alpha_{j;i}(p)$; this is similar to the universality found by Rudnick and Sarnak \cite{RS} in the $n$-level correlations. The higher moments of the Satake parameters control the rate of convergence to the random matrix predictions.} This leads to a very nice map from families of $L$-functions to $\{0,\pm 1\}$.\footnote{Instead of attaching a symmetry constant, additionally one can attach a symmetry vector which incorporates other information, such as the rank of the family.}

Another problem is that the main term in the $1$-level density agrees with random matrix theory, but the arithmetic of the family does not surface until we examine the lower order terms (which control the rate of convergence; see for example \cite{FI,Mil5,Yo1}). One promising line of research is the $L$-functions Ratios Conjecture \cite{CFZ1,CFZ2}, which is supported by corresponding calculations for random matrix ensembles (see \cite{CS,GJMMNPP,Mil4,Mil6,MilMo,Sto} for some recent work supporting these conjectures, especially \cite{CS} for a very accessible introduction to the method and a summary of its successes). Another approach is through hybrid product formulas \cite{GHK}. A typical $L$-function has two product representations, one as an Euler product over primes, and one as a Hadamard product over its zeros. In this approach an $L$-function is modeled by the product of a partial Euler and a partial Hadamard product. The Hadamard piece is believed to be well-modeled by random matrix theory, while the Euler product introduces the arithmetic. Thus the interplay between random matrix theory and number theory continues, and what began as a chance meeting in the 1970's now yields over 1,000,000 hits on a google search\footnote{The search string was 'number theory and random matrix theory'.} (as of August 2009).


\ \\

\end{document}